\documentclass{amsart}
\usepackage{amsmath,amssymb,amsfonts,verbatim}

\theoremstyle{plain}
\newtheorem{theorem}{Theorem}
\newtheorem{lemma}[theorem]{Lemma}
\newtheorem{proposition}[theorem]{Proposition}
\newtheorem{corollary}[theorem]{Corollary}
\def\A{\mathcal{E}}
\def\Ab{\bar{\rm A}}
\def\al{\alpha}
\def\Ar{{\rm A}}
\def\At{\widetilde A}

\def\be{\beta}
\def\bt{\ts\,\raise-0.5pt\hbox{\small$\boxtimes$}\,\,}
\def\Bt{\widetilde B}

\def\CC{\mathbb{C}}
\def\com{\ts,\hskip-.5pt}
\def\Cr{{\rm C}}
\def\Ct{\widetilde C}

\def\d{\partial}
\def\de{\delta}
\def\De{\Delta}
\def\Dt{\widetilde D}

\def\End{\operatorname{End}\ts}
\def\ep{\varepsilon}
\def\Ep{E^{\ts\prime}}

\def\ga{\gamma}
\def\ge{\geqslant}
\def\gl{\mathfrak{gl}}

\def\h{\mathfrak{h}}
\def\H{\mathfrak{H}}
\def\Hom{\operatorname{Hom}}

\def\I{{\rm I}}
\def\Ib{\tts\overline{\nns\rm I\nns}\tts}
\def\id{{\rm id}}
\def\io{\iota}

\def\J{{\rm J}}
\def\Jb{\,\overline{\!\rm J\ns}\ts}
\def\Jp{{\rm J}^{\ts\prime}}
\def\Jpb{\,\overline{\!\rm J\ns}\ts^{\ts\prime}}

\def\la{\lambda}
\def\lap{\la^{\ast}}
\def\lcd{\ts,\ldots,}
\def\le{\leqslant}

\def\mup{\mu^{\ast}}

\def\n{\mathfrak{n}}
\def\np{\n^{\ts\prime}}
\def\ns{\hskip-1pt}
\def\nns{\hskip-.5pt}

\def\om{\omega}

\def\op{\oplus}
\def\ot{\otimes}

\def\p{\mathfrak{p}}
\def\P{\mathcal{G}}
\def\Pb{\bar P}
\def\PD{\mathcal{GD}}

\def\q{\mathfrak{q}}
\def\Qb{\bar Q}
\def\qp{\q^{\ts\prime}}

\def\R{{\rm R}}
\def\Rb{\bar{\rm R}}

\def\S{\Lambda}
\def\si{\sigma}
\def\Sym{\mathfrak{S}}

\def\ts{\hskip1pt}
\def\tts{\hskip.5pt}

\def\U{\operatorname{U}}

\def\wp{w^{\ts\prime}}

\def\xic{\check\xi}

\def\Y{\operatorname{Y}}

\def\Z{\operatorname{Z}}
\def\ZZ{{\mathbb Z}} 

\begin{document} 
 
\title[Yangians and Mickelsson Algebras]{Yangians and Mickelsson Algebras II}

\author{Sergey Khoroshkin}
\address{
Institute for Theoretical and Experimental Physics, 
Moscow 117259,
Russia;}
\email{khor@itep.ru}

\author{Maxim Nazarov\\}
\address{
Department of Mathematics, 
University of York, 
York YO10 5DD, 
England;} 
\email{mln1@york.ac.uk}

\dedicatory{To Professor A.\,A.\,Kirillov
on the occasion of his 70\ts th birthday}



\begin{abstract}
We study a composition of two functors. The first one, 
from the category of modules over the Lie algebra $\gl_m$
to the category of modules over the degenerate affine Hecke algebra
of $GL_N\ts$, was introduced by I.\,Cherednik.
The second functor is a skew version
of the functor from the latter category to
the category of modules over the Yangian $\Y(\gl_n)$
due to V.\,Drinfeld.
We give a representation theoretic 
explanation of a link between intertwining operators on
tensor products of $\Y(\gl_n)$-modules, and
the \lq\lq\ts extremal cocycle\ts\rq\rq\  
on the Weyl group of $\gl_m$ introduced by D.\,Zhelobenko.
We also establish a connection between 
the composition of two functors, and
the \lq\lq\ts centralizer construction\ts\rq\rq\
of the Yangian $\Y(\gl_n)$ discovered by G.\,Olshanski.
\end{abstract}

\maketitle


\thispagestyle{empty} 


\renewcommand{\theequation}{\thesection.\arabic{equation}} 
\renewcommand{\thetheorem}{\thesection.\arabic{theorem}}
\makeatletter                                          
\@addtoreset{equation}{section}                        
\makeatother    


\thispagestyle{empty}

\section*{0.\ Introduction}

This article is a sequel to our work \cite{KN1} which
concerned two known functors. The definition of one of these
functors belongs to V.\,Drinfeld \cite{D2}. Let $\H_N$ be the
degenerate affine Hecke algebra corresponding to the
general linear group $GL_N$ over a non-Archimedean local field. This
is an associative algebra over the complex field $\CC$
which contains the symmetric group ring $\CC\,\Sym_N$ as a subalgebra.
Let $\Y(\gl_n)$ be the Yangian of the general linear Lie algebra $\gl_n\ts$.
This is a Hopf algebra over the field $\CC$
which contains the universal enveloping algebra $\U(\gl_n)$ as a subalgebra.
There is also a homomorphism of associative algebras 
$\pi_n:\Y(\gl_n)\to\U(\gl_n)$ identical on the subalgebra
$\U(\gl_n)\subset\Y(\gl_n)\ts$, see Section 1 of the present article
for the details. 
In \cite{D2} for any $\H_N$-module $W\ts$, an action of the algebra
$\Y(\gl_n)$ was defined on the vector space
$(W\ot(\CC^{\ts n})^{\ot N}){}^{\ts\Sym_N}$
of the diagonal $\Sym_N$-invariants in the tensor 
product of the vector spaces $W$ and $(\CC^{\ts n})^{\ot N}\ts$.
Thus one gets a functor from the category of all 
$\H_N$-modules to the category of $\Y(\gl_n)$-modules, called the
Drinfeld functor.

In \cite{KN1} we studied the composition of the Drinfeld functor
with another functor which was introduced by I.\,Cherednik \cite{C2}.
The latter functor was also studied by
T.\,Arakawa, T.\,Suzuki and A.\,Tsuchiya \cite{A,AS,AST}. 
For any module $V$ over the Lie algebra $\gl_m\ts$,
an action of the algebra $\H_N$ can be defined on the
tensor product of $\gl_m$-modules $V\ot(\CC^{\ts m})^{\ot N}\ts$. 
This action of $\H_N$
commutes with the diagonal action of $\gl_m$ on the tensor product.
Thus one gets a functor from the category of all 
$\gl_m$-modules to the category of $\H_N$-modules.
By applying the Drinfeld functor the $\H_N$-module 
$W=V\ot(\CC^{\ts m})^{\ot N}\ts$, one turns to an $\Y(\gl_n)$-module
the vector space
$$
(V\ot(\CC^{\ts m})^{\ot N}\ot(\CC^{\ts n})^{\ot N}){}^{\ts\Sym_N}
=
V\ns\ot\ts\operatorname{S\ts}^N(\CC^{\ts m}\ot\CC^{\ts n})\ts.
$$

In the present article we again use the Cherednik functor; 
details of its definition are reproduced in Section~1.
But we replace the Drinfeld functor by its skew version,
this version was also used in \cite{A}.
Similarly to \cite{D2}, for any $\H_N$-module $W$ an action of
the algebra $\Y(\gl_n)$ can be defined on the vector space
$(W\ot(\CC^{\ts n})^{\ot N}){}^{\ts\Sym_N}_{\,-}$
of the diagonal skew $\Sym_N$-invariants in the tensor 
product of $W$ and $(\CC^{\ts n})^{\ot N}\ts$.
Thus we obtain another functor from the category of all 
$\H_N$-modules to the category of $\Y(\gl_n)$-modules;
we call it the skew Drinfeld functor.
Details of its definition are also given in Section 1. 
By applying this functor
to $W=V\ot(\CC^{\ts m})^{\ot N}\ts$,
one turns to an $\Y(\gl_n)$-module the vector space
$$
(V\ot(\CC^{\ts m})^{\ot N}\ot(\CC^{\ts n})^{\ot N}){}^{\ts\Sym_N}_{\,-}
=
V\ns\ot\ts\S^N(\CC^{\ts m}\ot\CC^{\ts n})\ts.
$$

The action of 
$\Y(\gl_n)$ on this vector space
commutes with the action of 
$\gl_m\ts$. 
By taking the direct sum of these $\Y(\gl_n)$-modules
over $N=0,1\lcd m\ts n$
we turn into an $\Y(\gl_n)$-module
the space $V\ns\ot\ts\S\ts(\CC^{\ts m}\ot\CC^{\ts n})\ts$.
It is also a $\gl_m$-module; we denote this bimodule by
$\A_{\ts m}(V)\ts$. The additive group $\CC$ acts on 
the Hopf algebra $\Y(\gl_n)$ by automorphisms.
Let $\tau_z$ be the automorphism
of $\Y(\gl_n)$ corresponding to $z\in\CC\ts$, see \eqref{tauz}.
We denote by $\A_{\ts m}^{\,z}\ts(V)$ the $\Y(\gl_n)$-module obtained
from  $\A_{\ts m}(V)$ by pulling back through $\tau_{-z}\ts$.
By definition, $\A_{\ts m}^{\,z}\ts(V)$
coincides with $\A_{\ts m}(V)$ as a $\gl_{\ts m}\ts$-module.

Now take the Lie algebra $\gl_{m+l}\ts$.
Let $\p$ be the maximal parabolic subalgebra of $\gl_{m+l}$
containing the direct sum of Lie algebras $\gl_m\op\gl_{\ts l}\ts$.
Let $\q$ be the Abelian subalgebra of $\gl_{m+l}$ such that 
$
\gl_{m+l}=\q\op\p\ts.
$
For any $\gl_{\ts l}$-module $U$ let $V\bt U$ be the
$\gl_{m+l}$-module 
parabolically induced from the 
$\gl_m\op\gl_{\ts l}$-module $V\ot U\ts$.
This is a module induced from the subalgebra $\p\ts$.
Consider the space $\A_{\ts m+l}\ts(\ts V\ns\bt U\ts)_{\ts\q}$
of $\q$-coinvariants of the $\gl_{m+l}$-module 
$\A_{\ts m+l}\ts(\ts V\ns\bt U\ts)\ts$.
This space is an $\Y(\gl_n)$-module, which also
inherits the action of the Lie algebra $\gl_m\op\gl_{\ts l}\ts$.
Our Theorem \ref{parind} states that
the bimodule $\A_{\ts m+l}\ts(\ts V\ns\bt U\ts)_{\ts\q}$
over $\Y(\gl_n)$ and $\gl_m\op\gl_{\ts l}$
is equivalent to the tensor product 
$\A_{\ts m}(V)\ot\A_{\ts l}^{\ts m}\ts(U)\ts$.
Here we use the comultiplication on $\Y(\gl_n)\ts$.
This theorem is analogous to \cite[Theorem 2.1]{KN1}
which uses the original Drinfeld functor instead of its skew version.
Our proofs of the two theorems are also similar.

In Section 3 we establish 
a correspondence between intertwining operators 
on the tensor products of certain $\Y(\gl_n)$-modules, and
the \lq\lq\ts extremal cocycle\ts\rq\rq\ 
on the Weyl group $\Sym_m$ of the reductive Lie algebra $\ts\gl_m$ 
defined by D.\,Zhelobenko \cite{Z}.
Each of the tensor factors is obtained from one of the $\gl_n\ts$-modules 
$\S^N(\CC^{\ts n})$ where $N=1\com2\com\ts\ldots$
by pulling back through the homomorphism $\pi_n:\Y(\gl_n)\to\U(\gl_n)$
and then through one of the automorphisms $\tau_z$ of $\Y(\gl_n)\ts$.
In this article we assume that the parameters $z$ corresponding to
different tensor factors are in general position, that is their
differences do not belong to $\ZZ\ts$. It is well known that then
the tensor products are irreducible as $\Y(\gl_n)$-modules, see for
instance \cite{NT1}.

When each of the tensor factors
is obtained from one of the $\gl_n\ts$-modules 
$\operatorname{S\ts}^N(\CC^{\ts n})$
by pulling back through $\pi_n$
and then through one of the automorphisms $\tau_z\ts$,
such a correspondence was established
by V.\,Tarasov and A.\,Varchenko in \cite{TV2}. There
they used the classical duality theorem \cite{H}
which asserts that the images of the algebras
$\U(\gl_m)$ and $\U(\gl_n)$ in the ring 
of the differential operators on $\CC^{\ts m}\ot\CC^{\ts n}$
with polynomial coefficients are the commutants of each other.
Results relevant to this correspondence were also
obtained by Y.\,Smirnov and V.\,Tolstoy \cite{ST}.
In \cite{KN1} we gave a representation theoretic explanation
of this correspondence, by using the theory of
Mickelsson algebras \cite{M1,M2} as developed in \cite{KO}.
In the present article, we apply this theory to
the tensor products of the exterior powers of $\CC^{\ts n}$. 

We will identify the exterior algebra
$\S\ts(\ts\CC^{\ts m}\ot\CC^{\ts n})$
with the Grassmann algebra $\P(\ts\CC^{\ts m}\ot\CC^{\ts n})$
on $\ts m\ts n\ts$ anticommuting variables. The ring of endomorphisms
of the vector space $\P(\ts\CC^{\ts m}\ot\CC^{\ts n})$ is
denoted by $\PD\ts(\CC^{\ts m}\ot\CC^{\ts n})\ts$,
this ring is generated by all operators of left multiplication
by the anticommuting variables, and by the corresponding left derivations. 
For the details, see Section~1.
Now take the tensor product of associative algebras
\begin{equation}
\label{tenprod}
\U(\gl_m)\ot\PD\ts(\CC^{\ts m}\ot\CC^{\ts n})\ts.
\end{equation}
We have a representation
$\ga:\,\U(\gl_m)\to\PD\ts(\CC^{\ts m}\ot\CC^{\ts n})\ts$.
Taking the composition of the comultiplication map
on $\U(\gl_m)$ with the homomorphism $\id\ot\ga$
we get an embedding of $\U(\gl_m)$ into the algebra \eqref{tenprod}.
Our particular Mickelsson algebra is determined by the pair
formed by the algebra \eqref{tenprod}, and its subalgebra $\U(\gl_m)$
relative to this embedding. Other connections between the 
representation theory of Yangians and 
the theory of Mickelsson algebras were studied by A.\,Molev \cite{M}.

We complete this article with an observation on
the \lq\lq\ts centralizer construction\ts\rq\rq\ 
of the Yangian $\Y(\gl_n)$ due to G.\,Olshanski \cite{O1}.
For any two irreducible polynomial modules
$V$ and $V^{\ts\prime}$ of $\gl_m\ts$,
\cite{O1} provides an action
of $\Y(\gl_n)$ on the vector space
$$
\Hom_{\,\gl_m}(\ts V^{\ts\prime}\ts,
V\ot\S\ts(\ts\CC^{\ts m}\ot\CC^{\ts n}\ts))\ts.
$$
Moreover, this action is irreducible. Our Proposition \ref{arol}
states that the same action is inherited from the bimodule 
$\A_{\ts m}(V)=V\ns\ot\ts\S\ts(\CC^{\ts m}\ot\CC^{\ts n})$
over $\Y(\gl_n)$ and $\gl_m\ts$.

The first author was supported by the RFBR grant 05-01-01086 and 
by the joint JSPS-RFBR grant 05-01-02934.\ 
The second author was supported by the EPSRC grant C511166,
and by the EC grant MRTN-CT2003-505078.
This work was done when both authors stayed at the
Max Planck Institute for Mathematics in Bonn.
We are very grateful to the staff of the institute for their kind help and 
hospitality.


\section*{1.\ Skew Drinfeld functor}
\setcounter{section}{1}
\setcounter{equation}{0}
\setcounter{theorem}{0}

We begin with recalling two well known constructions
from the representation theory of the 
\textit{degenerate affine Hecke algebra\/} $\H_N\ts$, which
corresponds to the general linear group $GL_N$ over a local
non-Archimedean field. This algebra was introduced by V.\,Drinfeld
[D2], see also [L]. By definition, the complex associative algebra $\H_N$ is
generated by the symmetric group algebra $\CC\ts\Sym_N$ and by the pairwise
commuting elements $x_1\lcd x_N$ with the cross relations for $p=1\lcd N-1$
and $q=1\lcd N$
\begin{eqnarray}
\label{cross1}
\si_{p}\,x_q&=&x_q\,\si_{p}\ts,\quad q\neq p\ts,p+1\ts;
\\
\label{cross2}
\si_{p}\,x_p&=&x_{p+1}\,\si_{p}-1\ts.
\end{eqnarray}
Here and in what follows $\si_p\in \Sym_N$ denotes the 
transposition of numbers $p$ and $p+1\ts$. 
More generally, $\si_{pq}\in \Sym_N$ will denote the 
transposition of the numbers $p$ and $q\ts$. The group algebra
$\CC\ts\Sym_N$ can be then regarded as a subalgebra in $\H_N\ts$.
Furhtermore, it follows from the defining relations of $\H_N$
that a homomorphism $\H_N\to\CC\ts\Sym_N\ts$, identical
on the subalgebra $\CC\ts\Sym_N\subset \H_N\ts$,
can be defined by the assignments
\begin{equation}
\label{evalhecke}
x_p\mapsto\si_{1p}+\ldots+\si_{p-1,p}
\quad\text{for}\quad
p=1\lcd N.
\end{equation}

We will also use the elements of the algebra $\H_N\ts$,
$$
y_p=x_p-\si_{1p}-\ldots-\si_{p-1,p}
\quad\text{where}\quad
p=1\lcd N.
$$
Note that $y_p\mapsto0$ under the homomorphism $\H_N\to\CC\ts\Sym_N\ts$,
defined by
\eqref{evalhecke}. For any permutation $\si\in\Sym_N$, we have
\begin{equation}
\label{yrel}
\si\,y_p\,\si^{-1}=y_{\si(p)}\ts.
\end{equation}
It suffices to verify \eqref{yrel} when
$\si=\si_q$ and $q=1\lcd N-1$. Then (\ref{yrel})
is equivalent to the relations \eqref{cross1},\eqref{cross2}.
The elements $y_1\lcd y_N$ do not commute, but satisfy
the commutation relations
\begin{equation}
\label{ycom}
[y_p\ts,y_q]=\si_{pq}\cdot(y_p-y_q)\ts.
\end{equation}
Let us verify the equality in \eqref{ycom}.
Both sides of \eqref{ycom} are antisymmetric in $p$ and $q\ts$,
so it suffices to consider only the case when $p<q\ts$. Then
by using \eqref{yrel},
\begin{eqnarray}
\nonumber
&&[y_p\ts,x_q]=[x_p-\si_{1p}-\ldots-\si_{p-1,p}\ts,x_q]=0\ts,
\\
\nonumber
&&[y_p\ts,y_q]=[y_p\ts,y_q-x_q]=-\ts[y_p\ts,\si_{1q}+\ldots+\si_{q-1,q}]=
\si_{pq}\cdot(y_p-y_q)\,.
\end{eqnarray}
The algebra $\H_N$ is generated by $\CC\ts\Sym_N$
and the elements $y_1\lcd y_N\ts$.
Together with relations in $\CC\ts\Sym_N$, 
\eqref{yrel} and \eqref{ycom} are defining relations for $\H_N\ts$.
For more details on this presentation of the algebra $\H_N$ see 
for instance \cite[Section 1.3]{A}.

The first construction we recall here is due to I.\,Cherednik
\cite[Example 2.1]{C2}. It was further studied by
T.\,Arakawa, T.\,Suzuki and A.\,Tsuchiya \cite[Section~5.3]{AST}. 
Let $V$ be any module over the complex general linear Lie algebra $\gl_m\ts$.
Let $E_{ab}\in\gl_m$ with $a,b=1\lcd m$ be the
standard matrix units. We will also regard the matrix units $E_{ab}$
as elements of the algebra $\End(\CC^{\ts m})$, this should not cause any 
confusion. Let us consider the tensor product $V\ot(\CC^{\ts m})^{\ot N}$
of $\gl_m$-modules. Here each of the $N$ tensor factors $\CC^{\ts m}$ is
a copy of the natural $\gl_m$-module. We will use the indices
$1\lcd N$ to label these $N$ tensor factors.
For any $p=1\lcd N$
denote by $E^{\ts(p)}_{ab}$ the operator on the vector space 
$(\CC^{\ts m})^{\ot N}$ acting as
$$
\id^{\ts\ot\ts(p-1)}\ot E_{ab}\ot\id^{\ts\ot\ts(N-p)}\,.
$$
The following proposition coincides with \cite[Proposition 1.1]{KN1}.

\begin{proposition}
\label{ast}
{\rm\ (i)} 
An action of the algebra $\H_N$ on the vector space 
$V\ot(\CC^{\ts m})^{\ot N}$
can be defined as follows: 
the group $\Sym_N\subset \H_N$ acts naturally by
permutations of the $N$ tensor factors\/ $\CC^{\ts m}$, while any element 
$y_p\in \H_N$ acts as 
\begin{equation}
\label{yact}
\sum_{a,b=1}^m E_{\ts ba}\ot E_{ab}^{\ts(p)}.
\end{equation}
{\rm(ii)} 
This action of $\H_N$ commutes with the
{\rm(}diagonal\,{\rm)} action of\/ 
$\gl_m$ on $V\ot(\CC^{\ts m})^{\ot N}\ts$.
\end{proposition}

Let us now consider the \textit{triangular decomposition}
of the Lie algebra $\gl_m\ts$,
\begin{equation}
\label{tridec}
\gl_m=\n\op\h\op\np
\end{equation}
where $\h$ is the Cartan subalgebra of $\gl_m$ with the basis vectors
$E_{11}\lcd E_{mm}\ts$. Here $\n$ and $\np$ are the nilpotent
subalgebras spanned respectively by the elements $E_{\ts ba}$ and
$E_{ab}$ for all $a,b=1\lcd m$ such that $a<b\ts$.
For any $\gl_m$-module $W$, denote by
$W_{\n}$ the vector space
$W/\ts\n\cdot W$
of the coinvariants of the action of the subalgebra $\n\subset\gl_m$ on $W$.
The Cartan subalgebra $\h\subset\gl_m$ acts on the vector space $W_{\n}\ts$.

Now consider the tensor product 
$W=V\ot(\CC^{\ts m})^{\ot N}$ as a left module over the algebra $\H_N\ts$.
The action of $\H_N$ on this module commutes with the
action of the Lie algebra $\gl_m\ts$, 
and hence with the action of the subalgebra $\n\subset\gl_m\ts$.
So the space $(\ts V\ot(\CC^{\ts m})^{\ot N})_{\ts\n}$
of coinvariants of the action of $\n$ 
is a quotient of the $\H_N$-module $V\ot(\CC^{\ts m})^{\ot N}$.
Thus we have a functor from the category of all 
$\gl_m$-modules to the category of bimodules over $\h$ and $\H_N\ts$,
\begin{equation}
\label{zelefun}
V\mapsto(\ts V\ot(\CC^{\ts m})^{\ot N})_{\ts\n}\ts.
\end{equation}

Now take the \textit{Yangian}
$\Y(\gl_n)$ of the general linear Lie algebra $\gl_n\ts$.
The Yangian $\Y(\gl_n)$ is a deformation of the universal
enveloping algebra of the polynomial current Lie algebra $\gl_n[u]$
in the class of Hopf algebras, see for instance \cite{D1}.
The unital associative algebra $\Y(\gl_n)$ has a family of generators 
$$
T_{ij}^{\ts(1)},T_{ij}^{\ts(2)},\ts\ldots
\quad\text{where}\quad
i\ts,\ns j=1\lcd n\ts. 
$$
Defining relations for these generators
can be written in terms of the formal series
\begin{equation}
\label{tser}
T_{ij}(u)=
\de_{ij}+T_{ij}^{\ts(1)}u^{-\ns1}+T_{ij}^{\ts(2)}u^{-\ns2}+\,\ldots
\,\in\,\Y(\gl_n)\,[[\ts u^{-1}\ts]]\,.
\end{equation}
Here $u$ is the formal parameter. Let $v$ be another formal parameter.  
Then the defining relations in the associative algebra $\Y(\gl_n)$
can be written as
\begin{equation}
\label{yangrel}
(u-v)\cdot[\ts T_{ij}(u)\ts,T_{kl}(v)\ts]\ts=\;
T_{kj}(u)\ts T_{il}(v)-T_{kj}(v)\ts T_{il}(u)\,,
\end{equation}
where $i\com j\com k\com l=1\lcd n\ts$.
If $n=1$, then the algebra $\Y(\gl_n)$ is commutative.
The relations \eqref{yangrel} imply that for any $z\in\CC\,$, the assignments
\begin{equation}
\label{tauz}
\tau_z\ts:\,T_{ij}(u)\ts\mapsto\,T_{ij}(u-z)
\quad\textrm{for}\quad
i\com j=1\lcd n
\end{equation}
define an automorphism $\tau_z$ of the algebra $\Y(\gl_n)\ts$. 
Here each of the formal
power series $T_{ij}(u-z)$ in $(u-z)^{-1}$ should be re-expanded 
in $u^{-1}$, and the assignment \eqref{tauz} is a correspondence
between the respective coefficients of series in $u^{-1}$.

Now let $E_{ij}\in\gl_n$ with $i,j=1\lcd n$ be the
standard matrix units. We will also regard the matrix units $E_{ij}$
as elements of the algebra $\End(\CC^{\ts n})$, this should not cause any 
confusion. The Yangian $\Y(\gl_n)$ contains 
the universal enveloping algebra $\U(\gl_n)$ as a subalgebra\ts;
the embedding $\U(\gl_n)\to\Y(\gl_n)$ can be defined by the assignments
$$
E_{ij}\mapsto T_{ij}^{\ts(1)}
\quad\text{for}\quad
i\ts,\ns j=1\lcd n\ts. 
$$
There is a homomorphism $\pi_n:\Y(\gl_n)\to\U(\gl_n)$ identical
on the subalgebra $\U(\gl_n)\subset\Y(\gl_n)\ts$, 
it can be defined by the assignments
\begin{equation}
\label{pin}
\pi_n\ts:\,T_{ij}^{\ts(2)},T_{ij}^{\ts(3)},\ts\ldots
\,\mapsto\,0
\quad\text{for}\quad
i\ts,\ns j=1\lcd n\ts. 
\end{equation}
For further details on the definition of 
the algebra $\Y(\gl_n)$ see \cite[Chapter 1]{MNO}.

The second construction we recall here is a modification
of a construction due to V.\,Drinfeld \cite{D2},
which originally motivated his 
definition of the degenerate affine Hecke algebra $\H_N\ts$. 
For $p=1\lcd N$ let $E^{\ts(p)}_{ij}$ be the operator on the vector space 
$(\CC^{\ts n})^{\ot N}$ acting as
$$
\id^{\ts\ot\ts(p-1)}\ot E_{ij}\ot\id^{\ts\ot\ts(N-p)}\,.
$$
The group $\Sym_N$ acts on the tensor
product $(\CC^{\ts n})^{\ot N}$ from the left
by permutations of the $N$ tensor factors.
Let $W$ be any $\H_N$-module.
The group $\Sym_N$ also
acts from the left on $W$, 
via the embedding $\CC\ts\Sym_N\to \H_N$.
Consider the subspace
\begin{equation}
\label{sinv}
(W\ot(\CC^{\ts n})^{\ot N}){}^{\ts\Sym_N}_{\,-}
\subset 
W\ot(\CC^{\ts n})^{\ot N}
\end{equation}
of skew invariants with respect to the diagonal action of $\Sym_N\ts$.
On this subspace, each of the elements $\si_1\lcd\si_{N-1}$ of $\Sym_N$
acts as $-1$.
In the next proposition we use the convention that
$y_p^{\ts0}=1$, the identity element of the algebra $\CC\ts\Sym_N\ts$.

\begin{proposition}
\label{d}
One can define
an action of the algebra\/ $\Y(\gl_n)$ on the vector space
$(W\ot(\CC^{\ts n})^{\ot N}){}^{\ts\Sym_N}_{\,-}$ so that 
for any $s=0,1,2,\ts\ldots$ the generator
$T_{ij}^{\ts(s+1)}$ acts as
\begin{equation}
\label{tact}
\sum_{p=1}^N\,y_p^{\ts s}\ot E^{\ts(p)}_{ij}\ts.
\end{equation}
\end{proposition}

\begin{proof}
As an operator on the vector space $W\ot(\CC^{\ts n})^{\ot N}$,
\eqref{tact}
commutes with the diagonal action of $\Sym_N\ts$, due to the
relations \eqref{yrel} for the generators $y_1\lcd y_N$ of $\H_N\ts$.
So the restriction of the operator \eqref{tact} to the
subspace \eqref{sinv} is well defined. 

Following \cite{D2}, an action of the algebra
$\Y(\gl_n)$ can be defined on the subspace
\begin{equation}
\label{psinv}
(W\ot(\CC^{\ts n})^{\ot N}){}^{\ts\Sym_N}
\subset 
W\ot(\CC^{\ts n})^{\ot N}
\end{equation}
of invariants with respect to the diagonal action of $\Sym_N\ts$.
The generator $T_{ij}^{\ts(s+1)}$ of the algebra $\Y(\gl_n)$
acts on the subspace \eqref{psinv} as
$$
\sum_{p=1}^N\,(-\ts y_p)^{\ts s}\ot E^{\ts(p)}_{ij}\ts,
$$
see \cite[Proposition 1.2]{KN1}.
Now observe that the assignments
$$
\si_p\mapsto-\ts\si_p
\quad\text{and}\quad
x_q\mapsto-\ts x_q
$$
for all $p=1\lcd N-1$ and $q=1\lcd N$
define an automorphism of the algebra $\H_N\ts$,
see the defining relations \eqref{cross1},\eqref{cross2}.
Under this automorphism $y_q\mapsto-\ts y_q\ts$.
Let $W^{\ast}$ be the $\H_N$-module obtained
by pulling the action of $\H_N$ on $W$ back through this automorphism.
By substituting $W^{\ast}$ for $W$ in \eqref{psinv}
we get Proposition \ref{d}.
\end{proof}

\noindent\textit{Remark.}
When $s=0\ts$, the sum \eqref{tact} describes the action
of the element $E_{ij}\in\gl_n$ on the tensor product space 
$W\ot(\CC^{\ts n})^{\ot N}$, and hence on its subspace \eqref{sinv}.
Here each of the $N$ tensor factors $\CC^{\ts n}$ is regarded as a copy
of the natural $\gl_n$-module, and the action of $\gl_n$ on $W$ is trivial.
Hence the action of the Yangian $\Y(\gl_n)$ on the subspace
\eqref{sinv} as defined in Proposition \ref{d},
is compatible with the embedding $\U(\gl_n)\to\Y(\gl_n)\ts$.
\qed

\medskip\smallskip
Thus we obtain a functor from the category of all 
$\H_N$-modules to the category of $\Y(\gl_n)$-modules
\begin{equation}
\label{drinfun}
W\mapsto(W\ot(\CC^{\ts n})^{\ot N}){}^{\ts \Sym_N}_{\,-}\ts.
\end{equation}
We call it the \textit{skew Drinfeld functor} for the 
Yangian $\Y(\gl_n)\ts$.
Let us now apply this functor to the $\H_N$-module
$W=V\ot(\CC^{\ts m})^{\ot N}$ where $V$ is an arbitrary $\gl_m$-module;
see Proposition~\ref{ast}. The vector space of the resulting
$\Y(\gl_n)$-module is
$$
(V\ot(\CC^{\ts m})^{\ot N}\ot(\CC^{\ts n})^{\ot N}){}^{\ts\Sym_N}_{\,-}
=
V\ot((\CC^{\ts m}\ot\CC^{\ts n})^{\ot N}){}^{\ts\Sym_N}_{\,-}
$$
where the group $\Sym_N$ acts by permutations of
the $N$ tensor factors $\CC^{\ts m}\ot\CC^{\ts n}\ts$. Hence the
resulting vector space is 
\begin{equation}
\label{yangmod}
V\ns\ot\ts\S^N(\CC^{\ts m}\ot\CC^{\ts n})
\end{equation}
where we take the $N$-th exterior power of the vector space 
$\CC^{\ts m}\ot\CC^{\ts n}\ts$. Note that the Lie algebra $\gl_m$ also acts
on \eqref{yangmod} as the tensor product of two $\gl_m$-modules. 

We can identify the vector space $\CC^{\ts m}\ot\CC^{\ts n}$
with its dual, 
so that the standard basis vectors of $\CC^{\ts m}\ot\CC^{\ts n}$
are identified with the corresponding coordinate functions
$x_{ai}$ 
where
$a=1\lcd m$
and
$i=1\lcd n\ts$.
The exterior algebra $\S\ts(\CC^{\ts m}\ot\CC^{\ts n})$
is then identified with the \textit{Grassmann algebra}
$\P\ts(\CC^{\ts m}\ot\CC^{\ts n})$.
The latter algebra is generated by the elements $x_{ai}$
subject to the anticommutation relations 
$
x_{ai}\,x_{bj}=-\ts x_{bj}\,x_{ai}
$
for all indices $a\com b=1\lcd m$ and $i\com j=1\lcd n\ts$.
Let $\d_{ai}$ be the operator of left derivation on 
$\P\ts(\CC^{\ts m}\ot\CC^{\ts n})$
corresponding to the variable $x_{ai}\ts$,
this operator is also called the \textit{inner multiplication}
in $\P\ts(\CC^{\ts m}\ot\CC^{\ts n})$ corresponding to the
element $x_{ai}\ts$.

The ring of $\CC\ts$-endomorphisms of $\P\ts(\CC^{\ts m}\ot\CC^{\ts n})$ 
is generated by all operators of left multiplication by 
$x_{ai}$, and all operators $\d_{ai}\ts$;
see for instance \cite[Appendix~2.3]{H}.
This ring will be denoted by 
$\PD\ts(\CC^{\ts m}\ot\CC^{\ts n})$.
In this ring, we have the relations
\begin{equation}
\label{cliff}
x_{ai}\,\d_{bj}\ts+\,\d_{bj}\,x_{ai}\,=\,\de_{ab}\,\de_{ij}\,.
\end{equation}
Hence the ring $\PD\ts(\CC^{\ts m}\ot\CC^{\ts n})$
is isomorphic to the Clifford algebra
corresponding to the direct sum of the vector space
$\CC^{\ts m}\ot\CC^{\ts n}$ with its dual.
We can now describe the action of $\Y(\gl_n)$ on
the vector space \eqref{yangmod};
cf.\ \cite[Section 3]{A}.

\begin{proposition}
\label{dast}
{\rm\,(i)} 
For any $s=0,1,2,\ts\ldots$ the generator
$T_{ij}^{\ts(s+1)}$
acts on the $\Y(\gl_n)$-module \eqref{yangmod}
as the element of the tensor product\/ 
$\U(\gl_m)\ot\PD\ts(\CC^{\ts m}\ot\CC^{\ts n})\ts$,
\begin{equation}
\label{combact}
\sum_{c_0,c_1,\ldots,c_s=1}^m
E_{c_1c_0}\,E_{c_2c_1}\ldots\ts E_{c_sc_{s-1}}
\ot x_{c_0i}\,\d_{c_sj}\ts.
\end{equation}
When $s=0$, the first tensor factor in the summand in\/
\eqref{combact} is understood as $1$.
\\
{\rm\ (ii)} 
The action of $\Y(\gl_n)$ on \eqref{yangmod}
commutes with the {\rm(}diagonal\,{\rm)} action of\/ $\gl_m\ts$.
\end{proposition}

\begin{proof}
First consider the action of the sum \eqref{tact} on the
vector space $W\ot(\CC^{\ts n})^{\ot N}$ where 
$W=V\ot(\CC^{\ts m})^{\ot N}$\ts.
By substituting the sum \eqref{yact} for $y_p$ in \eqref{tact},
we then get the sum
\begin{equation}
\label{subst}
\sum_{p=1}^N\,
\Big(
\sum_{a,b=1}^m E_{\ts ba}\ot E_{ab}^{\ts(p)}
\,\Big){\!\!\ns\phantom{\big)}}^s\ot E^{\ts(p)}_{ij}
\end{equation}
acting on the vector space 
$V\ot(\CC^{\ts m})^{\ot N}\ot(\CC^{\ts n})^{\ot N}$.
Using the relations
$$
E_{ab}^{\ts(p)}\,E_{cd}^{\ts(p)}=\de_{bc}\,E_{ad}^{\ts(p)}
\quad\text{for}\quad
p=1\lcd N
$$
the sum \eqref{subst} can be rewritten as 
$$
\sum_{p=1}^N\ 
\sum_{c_0,c_1\ldots,c_s=1}^m\ 
E_{c_1c_0}
\ldots\ts 
E_{c_sc_{s-1}}\ot
E_{c_0c_s}^{\ts(p)}\ot
E^{\ts(p)}_{ij}\ts.
$$
To prove the part (i) of the proposition, 
it remains to observe that after identifying the subspace
\begin{equation}
\label{extn}
((\CC^{\ts m})^{\ot N}\ot(\CC^{\ts n})^{\ot N}){}^{\ts\Sym_N}_{\,-}
\subset
(\CC^{\ts m})^{\ot N}\ot(\CC^{\ts n})^{\ot N}
\end{equation}
with the subspace $\P^N(\CC^{\ts m}\ot\CC^{\ts n})$ of degree $N$
of in the Grassmann algebra, the operator
$$
\sum_{p=1}^N\,\,
E_{c_0c_s}^{\ts(p)}\ot
E^{\ts(p)}_{ij}
$$
on the subspace \eqref{extn}
gets identified with the operator $x_{c_0i}\,\d_{c_sj}$
on $\P^N(\CC^{\ts m}\ot\CC^{\ts n})\ts$.
Part (ii) of Proposition \ref{dast}
follows from the respective part of Proposition \ref{ast}.
\end{proof}

\noindent\textit{Remark.}
By definition, the basis element $E_{ab}\in\gl_m$ acts on the
space \eqref{yangmod}~as 
\begin{equation}
\label{eabact}
E_{ab}\ot1+
\sum_{k=1}^n\,
1\ot x_{ak}\,\d_{\ts bk}\ts.
\end{equation}
One can easily verify by direct calculation 
that the elements \eqref{combact} and \eqref{eabact}
of the algebra $\U(\gl_m)\ot\PD\ts(\CC^{\ts m}\ot\CC^{\ts n})$
commute with each other. Using \cite[Theorem~2]{H1},
one can show that the commutant
in the algebra $\U(\gl_m)\ot\PD\ts(\CC^{\ts m}\ot\CC^{\ts n})$
of all elements \eqref{eabact} with $a,b=1\lcd m$ 
is generated by the subalgebra $\Z(\gl_m)\ot1$ and
all elements of the form \eqref{combact};
cf.\ \cite[Section 2.1]{O2}.
Here $\Z(\gl_m)$ denotes
the centre of the universal enveloping algebra $\U(\gl_m)$.
This extends the fact
\cite[Section 4.2]{H} that the two families of operators
on the vector space $\P\ts(\CC^{\ts m}\ot\CC^{\ts n})\ts$,
\begin{equation}
\label{glmact}
\sum_{k=1}^n\,
x_{ak}\,\d_{\ts bk}
\quad\text{where}\quad
a,b=1\lcd m
\end{equation}
and
\begin{equation}
\label{glnact}
\sum_{c=1}^m\,
x_{ci}\,\d_{cj}
\quad\text{where}\quad
i,j=1\lcd n
\end{equation}
generate their mutual commutants in the algebra 
$\PD\ts(\CC^{\ts m}\ot\CC^{\ts n})$.
Here the operators \eqref{glmact} and \eqref{glnact}
describe the
actions on $\P\ts(\CC^{\ts m}\ot\CC^{\ts n})$ of the elements 
$E_{ab}\in\gl_m$ and $E_{ij}\in\gl_n$ respectively.
\qed

\medskip\smallskip
We will finish this section with an observation
on matrices with entries from the universal enveloping
algebra $\U(\gl_m)\ts$.
Let $E$ be the $m\times m$ matrix whose $ab\ts$-entry is 
the generator $E_{ab}\in\gl_m\ts$. Let $\Ep$ be the
transposed matrix.
Take the matrix inverse to $u-\Ep\ts$.
Here the summand $u$ stands for the scalar $m\times m$ matrix
with diagonal entry $u\ts$, and the inverse is a formal power
series in $u^{-1}$ with matrix coefficients. 
Denote by $X_{ab}(u)$ the $ab\ts$-entry of inverse matrix. Then
$$
X_{ab}(u)=
\ts\de_{ab}\,u^{-1}\ns+
E_{\ts ba}\ts u^{-2}\,+\,
\sum_{s=1}^\infty\,
\sum_{c_1,\ldots,c_s=1}^m
E_{c_1a}\ts E_{c_2c_1}\ldots\ts E_{c_sc_{s-1}}\ts E_{\ts bc_s}
u^{-s-2}\ts.
$$
The assignment of the element \eqref{combact}
to any coefficient $T_{ij}^{\ts(s+1)}$
of the series \eqref{tser} can be now written as
$$
T_{ij}(u)\mapsto\de_{ij}\ts+\sum_{a,b=1}^m
X_{ab}(u)\ot x_{ai}\,\d_{\ts bj}\ts.
$$


\section*{2.\ Parabolic induction}
\setcounter{section}{2}
\setcounter{equation}{0}
\setcounter{theorem}{0}

The Yangian $\Y(\gl_n)$ is a Hopf algebra over the field $\CC\ts$.
Using the series \eqref{tser},
the comultiplication $\De:\Y(\gl_n)\to\Y(\gl_n)\ot\Y(\gl_n)$ is defined by
the assignment
\begin{equation}\label{1.33}
\De:T_{ij}(u)\ts\mapsto\ts\sum_{k=1}^n\ T_{ik}(u)\ot T_{kj}(u)\,;
\end{equation}
the tensor product at the right hand side of the assignment (\ref{1.33})
is taken over the subalgebra 
$\CC[[u^{-1}]]\subset\Y(\gl_n)\,[[u^{-1}]]\ts$.
When taking tensor products of modules over $\Y(\gl_n)$, 
we will use the comultiplication \eqref{1.33}.
The counit homomorphism 
$\ep:\Y(\gl_n)\to\CC$ is defined by 
$$
\ep:\,T_{ij}(u)\ts\mapsto\ts\de_{ij}\cdot1\ts.
$$
The antipode ${\rm S}$ on $\Y(\gl_n)$ is defined by using the 
$n\times n$ matrix $T(u)$ whose $ij$-entry is the series $T_{ij}(u)\ts$.
This matrix is invertible as formal power series in $u^{-1}$ with
matrix coefficients, because the leading term of this series is
the identity $n\times n$ matrix. Then
the involutive anti-automorphism ${\rm S}$ of $\Y(\gl_n)$
is defined by the assignment
$$
{\rm S}\ts:\ts T(u)\mapsto T(u)^{-1}.
$$
This assignment means that by applying ${\rm S}$ to the coefficients of
the series $T_{ij}(u)$, we obtain the series which is the
$ij$-entry of the inverse matrix $T(u)^{-1}\ts$. 
We also use the involutive automorphism $\om_n$
of the algebra $\Y(\gl_n)$ defined by a similar assignment,
\begin{equation}
\label{1.51}
\om_n:\ts T(u)\mapsto T(-u)^{-1}.
\end{equation}
For more details on the Hopf algebra
structure on $\Y(\gl_n)$ see \cite[Chapter 1]{MNO}.

Let us now consider the direct sum of bimodules over
$\gl_m$ and $\Y(\gl_n)\ts$,
$$ 
\mathop{\op}\limits_{N=0}^{mn}\,V\ns\ot\ts\S^N(\CC^{\ts m}\ot\CC^{\ts n})
=V\ns\ot\ts\S\,(\CC^{\ts m}\ot\CC^{\ts n})\ts.
$$
Let us denote this bimodule by $\A_{\ts m}(V)\ts$, so that $\A_{\ts m}$
is a functor from the category of all $\gl_m$-modules to the
category of bimodules over $\gl_m$ and $\Y(\gl_n)\ts$. 
By identifying the exterior algebra
$\S\,(\CC^{\ts m}\ot\CC^{\ts n})$ with the Grassmann algebra 
$\P\ts(\CC^{\ts m}\ot\CC^{\ts n})$, 
the action of the generator $T_{ij}^{\ts(s+1)}$ of $\Y(\gl_n)$
on $\A_{\ts m}(V)$ is described by \eqref{combact}.

For any positive integer $l$ let $U$ be a
module over the Lie algebra $\gl_{\ts l}\ts$. Then $\A_{\ts l}(U)$ 
is another $\Y(\gl_n)$-module. For any $z\in\CC$ let us denote by
$\A_{\ts l}^{\ts z}\ts(U)$ the $\Y(\gl_n)$-module obtained
from  $\A_{\ts l}(U)$ by pulling back through the automorphism
$\tau_{-z}$ of $\Y(\gl_n)\ts$, see \eqref{tauz}. 
As a $\gl_{\ts l}\ts$-module $\A_{\ts l}^{\ts z}\ts(U)$
coincides with $\A_{\ts l}(U)\ts$.

The decomposition $\CC^{\ts m+l}=\CC^{\ts m}\op\CC^{\ts l}$ 
determines an embedding
of the direct sum $\gl_m\op\gl_{\ts l}$ of Lie algebras into $\gl_{m+l}\ts$.
As a subalgebra of $\gl_{m+l}\ts$,
the direct summand $\gl_m$ is spanned by the matrix units 
$E_{ab}\in\gl_{m+l}$ where
$a,b=1\lcd m\ts$. The direct summand $\gl_{\ts l}$ is spanned by
the matrix units $E_{ab}$ where $a,b=m+1\lcd m+l\ts$.
Let $\q$ and $\qp$ be 
the Abelian subalgebras of $\gl_{m+l}$
spanned respectively by matrix units $E_{\ts ba}$ and 
$E_{ab}$ for all $a=1\lcd m$ and $b=m+1\lcd m+l\ts$.
Put $\p=\gl_m\op\gl_{\ts l}\op\qp\ts$.
Then $\p$ is a maximal parabolic subalgebra of the reductive 
Lie algebra $\gl_{m+l}\ts$, and moreover 
$
\gl_{m+l}=\q\op\p\ts.
$
Denote by $V\ns\bt U$ the $\gl_{m+l}$-module \textit{parabolically induced\/}
from the $\gl_m\op\gl_{\ts l\ts}$-module $V\ot U$.
To define $V\ns\bt U$, one first extends the action of the Lie algebra
$\gl_m\op\gl_{\ts l\ts}$ on $V\ot U$ to the Lie algebra $\p\ts$, so that
any element of the subalgebra $\qp\subset\p$ acts on $V\ot U$ as zero.
By definition, $V\ns\bt U$ is the $\gl_{m+l}$-module induced from the 
$\p$-module $V\ot U$.  

Let us now consider the bimodule $\A_{\ts m+l}\ts(\ts V\ns\bt U\ts)$
over $\gl_{m+l}$ and $\Y(\gl_n)\ts$.
Here the action of $\Y(\gl_n)$ commutes with the
action of the Lie algebra $\gl_{m+l}\ts$, 
and hence with the action of the subalgebra $\q\subset\gl_{m+l}\ts$.
Therefore the vector space 
$\A_{\ts m+l}\ts(\ts V\ns\bt U\ts)_{\ts\q}$
of coinvariants of the action of the subalgebra $\q$ 
is a quotient of the $\Y(\gl_n)$-module 
$\A_{\ts m+l}\ts(\ts V\ns\bt U\ts)\ts$.
The subalgebra 
$\gl_m\op\gl_{\ts l}\subset\gl_{m+l}$ also acts on this quotient space.

\begin{theorem}
\label{parind}
The bimodule $\A_{\ts m+l}\ts(\ts V\ns\bt U\ts)_{\ts\q}$
over the Yangian $\Y(\gl_n)$ and
the direct sum $\gl_m\op\gl_{\ts l}\ts$, 
is equivalent to the tensor product 
$\A_{\ts m}(V)\ot\A_{\ts l}^{\ts m}\ts(U)\ts$.
\end{theorem}

Our proof of the theorem is based on two simple lemmas.
The first of them
applies to matrices over arbitrary unital ring.
Take a matrix of size $(m+l)\times(m+l)$ over such a ring,  
and write it as the block matrix
\begin{equation}
\label{blockmat}
\begin{bmatrix}\,A\,&B\,\\\,C\,&D\,\end{bmatrix}
\end{equation}
where the blocks $A,B,C,D$ are matrices of sizes
$m\times m$, $m\times l$, $l\times m$, $l\times l$ respectively. 
The following fact is well known, see for instance \cite[Lemma 2.2]{KN1}.

\begin{lemma}
\label{lemma1}
Suppose that the matrix \eqref{blockmat} is invertible.
Suppose that the matrices $A$ and $D$ are also invertible. Then the matrices
$A-B\ts D^{-1}\ts C$ and $D-C\ts A^{-1}B$ are invertible too, and
$$
\begin{bmatrix}A&B\\C&D\end{bmatrix}^{-1}\!=\ \  
\begin{bmatrix}
(A-B\ts D^{-1}\ts C)^{-1}
&
\,-A^{-1}B\ts(\ts D-C\ts A^{-1}B)^{-1}\,
\\
\,-D^{-1}C\ts(A-B\ts D^{-1}\ts C)^{-1}\,
&
(\ts D-C\ts A^{-1}B)^{-1}
\end{bmatrix}
\,.
$$
\end{lemma}

Consider again
the $m\times m$ matrix $E$ whose $ab\ts$-entry is 
the generator $E_{ab}\in\gl_m\ts$. The $ab\ts$-entry 
of the matrix inverse to $u-\Ep$ has been denoted by
$X_{ab}(u)\ts$.
Denote by $Z(u)$ the trace of the inverse matrix, so that
\begin{equation}
\label{zu}
Z(u)\,=\,\sum_{c=1}^m\,X_{cc}(u)\ts.
\end{equation}
Then $Z(u)$ is a formal power series in $u^{-1}$
with the coefficients from the algebra  $\U(\gl_m)\ts$.
Note that the leading term of this series is $m\ts u^{-1}$.
Let us now regard the coefficents of the series $X_{ab}(u)$ and
$Z(u)$ as elements of the algebra
$\U(\gl_{m+l})\ts$, using the standard embedding of 
the Lie algebra $\gl_m$ to $\gl_{m+l}\ts$.
Our next lemma can be obtained by replacing 
in \cite[Lemma 2.3]{KN1} the formal parameter $u$ by $-\ts u\ts$.

\begin{lemma}
\label{lemma2}
For any\/ $a=1\lcd m$ and\/ $d=1\lcd l$ 
we have an equality of the series
with coefficients in the algebra\/ $\U(\gl_{m+l})\ts,$
\begin{equation}
\label{xe}
\sum_{b=1}^m\,E_{\ts m+d,b}\,X_{ab}(u)
\,=\,
\sum_{b=1}^m\,X_{ab}(u)\,E_{\ts m+d,b}\,(1+Z(u))\ts.
\end{equation}
\end{lemma}

\smallskip\medskip
\noindent\textit{Proof of Theorem \ref{parind}.}
The 
vector space of $\gl_{m+l}\ts$-module $V\bt U$
can be identified with the tensor product $\U(\q)\ot V\ns\ot U$
so that the Lie subalgebra $\q\subset\gl_{m+l}$ acts via
left multiplication on the first tensor factor.
Note that the corresponding
action of the commutative algebra $\U(\q)$ is free.
The tensor product $V\ns\ot U$ is then identified with the subspace
\begin{equation}
\label{onesub}
1\ot V\ns\ot U\subset \U(\q)\ot V\ns\ot U\ts.
\end{equation}
On this subspace,
any element of the subalgebra $\qp\subset\gl_{m+l}$ acts as zero, 
while the two direct summands
of subalgebra $\gl_m\op\gl_{\ts l}\subset\gl_{m+l}$
act non-trivially only on the tensor factors $V$ and $U$ respectively.
All this determines the action of Lie algebra $\gl_{m+l}$ on
$\U(\q)\ot V\ns\ot U$. Now consider $\A_{\ts m+l}\ts(\ts V\ns\bt U\ts)$
as a $\gl_{m+l}\ts$-module, we will denote it by $W$ for short.
Then $W$ is the tensor product of two $\gl_{m+l}\ts$-modules,
$$
W=(V\bt U)\ot\P\ts(\CC^{\ts m+l}\ot\CC^{\ts n})=
\U(\q)\ot V\ns\ot U\ot\P\ts(\CC^{\ts m+l}\ot\CC^{\ts n})\ts.
$$

The vector spaces of the two $\Y(\gl_n)$-modules 
$\A_{\ts m}(V)$ and $\A_{\ts l}^{\ts m}\ts(U)$ are respectively
$$
V\ns\ot\ts\P\,(\CC^{\ts m}\ot\CC^{\ts n})
\ \quad\text{and}\ \quad
U\ns\ot\ts\P\,(\CC^{\ts l}\ot\CC^{\ts n})\ts.
$$
Identify the tensor product of these two vector spaces with
\begin{equation}
\label{vuprod}
V\ns\ot U\ot
\P\,(\CC^{\ts m}\ot\CC^{\ts n})\ot\P\,(\CC^{\ts l}\ot\CC^{\ts n})=
V\ns\ot U\ot
\P\,(\CC^{\ts m+l}\ot\CC^{\ts n})
\end{equation}
where we use the standard direct sum decomposition 
$$\CC^{\ts m+l}\ot\CC^{\ts n}=
\CC^{\ts m}\ot\CC^{\ts n}\op\ts\CC^{\ts l}\ot\CC^{\ts n}\ts.
$$
Regard the tensor product $V\ns\ot U$ in \eqref{vuprod} as a module of
the subalgebra $\gl_m\op\gl_{\ts l}\subset\gl_{m+l}$\ts. This subalgebra
also acts on $\P\,(\CC^{\ts m+l}\ot\CC^{\ts n})$ naturally.
Define a linear map
$$
\chi:\,V\ns\ot U\ot\P\,(\CC^{\ts m+l}\ot\CC^{\ts n})\ts\to W\ts/\,\q\cdot W
$$
by the assignment
$$
\chi:\,y\ot x\ot f\,\mapsto\,1\ot y\ot x\ot f\,+\,\q\cdot W
$$
for any $y\in V$, $x\in U$ and $f\in\P\,(\CC^{\ts m+l}\ot\CC^{\ts n})\ts$. 
The operator $\chi$ evidently intertwines the actions of the Lie
algebra $\gl_m\op\gl_{\ts l}\ts$. 


Let us demonstrate that the operator $\chi$ is bijective.
Firstly consider the action of the Lie subalgebra $\q\subset\gl_{m+l}$
on the vector space
$$
\P\ts(\CC^{\ts m+l})\ts=\ts
\P\ts(\CC^{\ts m})\ot\P\ts(\CC^{\ts l})\ts.
$$
This vector space admits an ascending filtration by the subspaces
$$
\mathop{\op}\limits_{N=0}^{\ts K}  
\P^{\ts N}(\CC^{\ts m})\ot\P\,(\CC^{\ts l})
\ \quad\textrm{where}\ \quad
K=0,1\lcd m\ts.
$$
Here $\P^{\ts N}(\CC^{\ts m})$ is
the subspace of degree $N$ in the Grassmann algebra $\P\ts(\CC^{\ts m})\ts$.
The action of the Lie algebra $\q$ on $\P\ts(\CC^{\ts m+l})$
preserves each of these subspaces, and is trivial
on the associated graded space. Similarly, the vector
space $\P\,(\CC^{\ts m+l}\ot\CC^{\ts n})$ admits an ascending
filtration by $\q$-submodules such that $\q$ acts trivially on
each of the corresponding graded subspaces.
The latter filtration induces a filtration of $W$
by $\q$-submodules such that the corresponding graded quotient
$\operatorname{gr}W$ is a free $\U(\q)$-module. The
space of coinvariants $(\ts\operatorname{gr}W)_{\ts\q}$ is therefore
isomorphic to the vector space
$V\ns\ot U\ot\P\,(\CC^{\ts m+l}\ot\CC^{\ts n})\ts$,
via the bijective linear map
$$
y\ot x\ot f\,\mapsto\,1\ot y\ot x\ot f\,+\,\q\cdot(\operatorname{gr}W)\ts.
$$
Therefore the linear map $\chi$ is bijective as well.

Let us now prove that the map $\chi$ intertwines the actions of
the Yangian $\Y(\gl_n)\ts$. Take the $(m+l)\times(m+l)$ matrix whose
$ab\ts$-entry is $\de_{ab}\ts u-\ns E_{\ts ba}\ts$; we regard 
$E_{ba}$ as an element of $\U(\gl_{m+l})\ts$. 
Write this matrix in the block form 
\eqref{blockmat} where $A,B,C,D$ are matrices of sizes
$m\times m$, $m\times l$, $l\times m$, $l\times l$ respectively.
In the notation introduced in the end of Section 1, 
here $A=u-\ns\Ep$.  
Using the observation made there along with the definition \eqref{1.33}
of the comultiplication, the action of the algebra $\Y(\gl_N)$ on
the vector space \eqref{vuprod} of the tensor product of the
$\Y(\gl_n)$-modules $\A_{\ts m}(V)$ and $\A_{\ts l}^{\ts m}\ts(U)$
can be described by assigning to every series $T_{ij}(u)$ the 
product of the series 
$$
\sum_{k=1}^n\,\,
\Bigl(\ts\de_{ik}\,+\!
\sum_{a,b=1}^m
(A^{-1})_{ab}\ot x_{ai}\,\d_{\ts bk}\ts\Bigr)\,\times
$$
$$
\Bigl(\ts\de_{kj}\,+\!
\sum_{c,d=1}^{l}
((\ts D+m\ts )^{-1})_{cd}\ot x_{\ts m+c,k}\,\d_{\ts m+d,j}\ts\Bigr)\,=
$$
\begin{equation}
\label{sum1}
\de_{ij}\,+\!\ts
\sum_{a,b=1}^m
(A^{-1})_{ab}\ot x_{ai}\,\d_{\ts bj}\,+
\sum_{c,d=1}^{l}
((\ts D+m\ts )^{-1})_{cd}\ot x_{\ts m+c,i}\,\d_{\ts m+d,j}\ +
\end{equation}
\begin{equation}
\label{sum2}
\sum_{k=1}^n\ \,
\sum_{a,b=1}^m\,
\sum_{c,d=1}^{l}\ \,
(A^{-1})_{ab}\,((\ts D+m\ts )^{-1})_{cd}\ot 
x_{ai}\,\d_{\ts bk}\,x_{\ts m+c,k}\,\d_{\ts m+d,j}\,.
\end{equation}
Note that in \eqref{sum2} we have
$\d_{\ts bk}\,x_{\ts m+c,k}=-\ts x_{\ts m+c,k}\,\d_{\ts bk}$ 
because $b\le m\ts$; see \eqref{cliff}.
The first tensor factors of all summands in
\eqref{sum1} and \eqref{sum2} correspond to the
action of the universal enveloping algebra
$\U(\ts\gl_m\op\gl_{\ts l}\ts)$ on $V\ns\ot U$.

Let us now write the matrix inverse to \eqref{blockmat} as
the block matrix
$$
\begin{bmatrix}
\ \At\,&\Bt\ 
\\
\ \Ct\,&\Dt\ 
\end{bmatrix}
$$
where $\At,\Bt,\Ct,\Dt$ are matrices of sizes
$m\times m$, $m\times l$, $l\times m$, $l\times l$ respectively.
Each of these four blocks is regarded 
as formal power series in $u^{-1}$ 
with matrix coefficients. The entries of these matrix coefficients 
belong to $\U(\ts\gl_{m+l})\ts$.
By once again using the observation made in the end
of Section 1, the action of $\Y(\gl_n)$
on the vector space $W$ can now be described by assigning to
every series $T_{ij}(u)$ the sum
$$
\de_{ij}+\sum_{a,b=1}^m
\At_{\ts ab}\ot x_{ai}\,\d_{\ts bj}\,+\,
\sum_{a=1}^m\,\sum_{d=1}^{l}\, 
\Bt_{\ts ad}\ot x_{ai}\,\d_{\ts m+d,j}\ +
$$
$$
\sum_{b=1}^m\,\sum_{c=1}^{l}\, 
\Ct_{\ts cb}\ot x_{\ts m+c,i}\,\d_{\ts bj}\,+
\sum_{c,d=1}^{l}
\Dt_{\ts cd}\ot x_{\ts m+c,i}\,\d_{\ts m+d,j}\,.
$$
The first tensor factors in the summands correspond to the
action of the algebra $\U(\gl_{m+l})$ on the space
$\U(\q)\ot V\ns\ot U$ of the parabolically induced module $V\bt U\ts$.

Apply these tensor factors to elements of the
subspace \eqref{onesub}. By Lemma \ref{lemma1}, 
$$
\At=(A-B\ts D^{-1}\ts C)^{-1}\ts.
$$
All entries of the matrix $C$ belong to $\qp$ and
hence act on the subspace \eqref{onesub} as zeroes.
Further, we have $A=u-\Ep$. Every entry of the matrix $\Ep$
belongs to the subalgebra $\gl_m\subset\gl_{m+l}$ and
the adjoint action of this subalgebra on $\gl_{m+l}$
preserves $\qp$. Therefore the results of applying
$(A^{-1})_{ab}$ and $\At_{ab}$ to elements of
the subspace \eqref{onesub} are the same.
Similar arguments show that any entry of the matrix
$$
\Ct=-\ts D^{-1}C\ts(A-B\ts D^{-1}\ts C)^{-1}=-\ts D^{-1}C\ts\At
$$
act on the subspace \eqref{onesub} as zero.

Consider the matrix 
$$
\Dt=(\ts D-C\ts A^{-1}B)^{-1}\ts.
$$
In the notation of Lemma \ref{lemma2} the $ab\ts$-entry of
$A^{-1}$ is $X_{ab}(u)\ts$, and the trace of $A^{-1}$ is $Z(u)\ts$. 
Using that lemma, the $cd\ts$-entry of the $l\times l$ 
matrix $D-C\ts A^{-1}B$ equals
$$
\de_{cd}\,u\ts-E_{m+d,m+c}-\sum_{a,b=1}^m 
E_{\ts a,m+c}\,X_{ab}(u)\,E_{\ts m+d,b}\,=
\de_{cd}\,u\ts-E_{m+d,m+c}
$$
$$
-\sum_{a,b=1}^m 
E_{\ts a,m+c}\,E_{\ts m+d,b}\,X_{ab}(u)\,(\ts1+Z(u))^{-1}=\,
\de_{cd}\,u\ts-E_{m+d,m+c}
\vspace{4pt}
$$
$$
-\sum_{a,b=1}^m 
(\ts E_{\ts m+d,b}\,E_{\ts a,m+c}+
\de_{cd}\,E_{ab}-\de_{ab}\,E_{\ts m+d,m+c}\ts)
\,X_{ab}(u)\,(\ts1+Z(u))^{-1}=
\vspace{7pt}
$$
\begin{equation}
\label{dmcd}
\ \de_{cd}\,(\ts u+(\ts m-u\,Z(u))\ts(\ts1+Z(u))^{-1})+
E_{m+d,m+c}\,(\ts Z(u)\ts(\ts1+Z(u))^{-1}-\ts1\ts)
\vspace{8pt}
\end{equation}
\begin{equation}
\label{efactor}
-\,
\sum_{a,b=1}^m 
E_{\ts m+d,b}\,E_{\ts a,m+c}
\,X_{ab}(u)\,(\ts1+Z(u))^{-1}.
\vspace{10pt}
\end{equation}
We used the identity
$$
\sum_{a,b=1}^m E_{ab}\,X_{ab}(u)\,=\,-\,m+u\,Z(u)
$$
which follows from \eqref{zu}. The expression \eqref{dmcd} equals
$$
(\ts D+m)_{\ts cd}\ts(\ts1+Z(u))^{-1}.
$$
The factor $E_{\ts a,m+c}$ in any summand in \eqref{efactor}
belongs to $\qp$ while every element of $\qp$
acts on the subspace \eqref{onesub} as zero. 
The coefficients of the series $X_{ab}(u)$ and $Z(u)$ in \eqref{efactor}
belong to $\U(\gl_m)$ while the adjoint action of subalgebra 
$\gl_m\subset\gl_{m+l}$ preserves $\qp$. The adjoint action of
the element $E_{m+d,m+c}\in\gl_{m+l}$ also preserves $\qp$.
Hence the result of applying $\Dt_{cd}$ to elements of
the subspace \eqref{onesub} is the same as that of applying
$$
(\ts1+Z(u))\ts((\ts D+m)^{-1})_{cd}\ts.
$$

Now consider
$$
\Bt=-\ts A^{-1}B\ts(\ts D-C\ts A^{-1}B)^{-1}=-\ts A^{-1}B\ts\Dt\ts.
$$
The above arguments show that the result of applying
the $ad\ts$-entry of this matrix to
elements of the subspace \eqref{onesub} is the same as that of applying
the $ad\ts$-entry of
$$
-\ts A^{-1}B\ts(\ts1+Z(u))\ts(\ts D+m)^{-1}\ts.
$$
By using Lemma \ref{lemma2} once again, the latter entry equals
$$
\sum_{b=1}^m\,\sum_{c=1}^l\,
X_{ab}(u)\,E_{\ts m+c,b}
(\ts1+Z(u))\ts((\ts D+m)^{-1})_{cd}\,=
$$
$$
\sum_{b=1}^m\,\sum_{c=1}^l\,
E_{\ts m+c,b}\,X_{ab}(u)\ts((\ts D+m)^{-1})_{cd}\,.
\vspace{4pt}
$$

Thus we have proved that
the action of $\Y(\gl_n)$ on the elements of the subspace
\begin{equation}
\label{1vu}
1\ot V\ns\ot U\ot\P\ts(\CC^{\ts m+l}\ot\CC^{\ts n})\subset W
\end{equation}
can be described by assigning to every series $T_{ij}(u)$ 
the sum of the series 
$$
\de_{ij}+\sum_{a,b=1}^m
(A^{-1})_{\ts ab}\ot x_{ai}\,\d_{\ts bj}\,+
\sum_{c,d=1}^{l}
(1+Z(u))\ts((\ts D+m)^{-1})_{\ts cd}\ot x_{\ts m+c,i}\,\d_{\ts m+d,j}\ +
\vspace{-6pt}
$$
\begin{equation}
\label{exud}
\sum_{a,b=1}^m\,\sum_{c,d=1}^{l}\, 
E_{\ts m+c,b}\,X_{ab}(u)\ts((\ts D+m)^{-1})_{cd}
\ot x_{ai}\,\d_{\ts m+d,j}\ts.
\vspace{4pt}
\end{equation}
Let us now consider the results of the action of $\Y(\gl_n)$
on this subspace modulo $\q\cdot W$. Since $E_{\ts m+c,b}\in\q\,$,
the expession
in the line \eqref{exud} can be then replaced by
$$
-\ \sum_{k=1}^n\ \sum_{a,b=1}^m\,\sum_{c,d=1}^{l}\, 
X_{ab}(u)\ts((\ts D+m)^{-1})_{cd}
\ot x_{m+c,k}\,\d_{\ts bk}\,x_{ai}\,\d_{\ts m+d,j}\,=
$$
$$
\sum_{k=1}^n\ \sum_{a,b=1}^m\,\sum_{c,d=1}^{l}\, 
(A^{-1})_{ab}\ts((\ts D+m)^{-1})_{cd}
\ot x_{m+c,k}\,x_{ai}\,\d_{\ts bk}\,\d_{\ts m+d,j}
$$
$$
-\,\sum_{c,d=1}^{l}\, 
Z(u)\ts((\ts D-m)^{-1})_{cd}
\ot x_{m+c,i}\,\d_{\ts m+d,j}\ts.
$$
Here we used the equality of differential operators
$\d_{\ts bk}\,x_{ai}=
-x_{ai}\,\d_{\ts bk}+\de_{ab}\,\de_{ik}\ts$;
see \eqref{cliff}.
By making this replacement we show 
that modulo $\q\cdot W$,
the action of $\Y(\gl_n)$ on elements of
the subspace \eqref{1vu}
can be described by assigning to every series $T_{ij}(u)$ 
the sum of the series 
$$
\de_{ij}+\sum_{a,b=1}^m
(A^{-1})_{\ts ab}\ot x_{ai}\,\d_{\ts bj}\,+
\sum_{c,d=1}^{l}
(1+Z(u))\ts((\ts D+m)^{-1})_{\ts cd}\ot x_{\ts m+c,i}\,\d_{\ts m+d,j}\ +
\vspace{-2pt}
$$
$$
\sum_{k=1}^n\ \sum_{a,b=1}^m\,\sum_{c,d=1}^{l}\, 
(A^{-1})_{ab}\ts((\ts D+m)^{-1})_{cd}
\ot x_{m+c,k}\,x_{ai}\,\d_{\ts bk}\,\d_{\ts m+d,j}
$$
$$
-\,\sum_{c,d=1}^{l}\, 
Z(u)\ts((\ts D+m)^{-1})_{cd}
\ot x_{m+c,i}\,\d_{\ts m+d,j}\,.
$$
Using the equality $x_{m+c,k}\,x_{ai}=-\ts x_{ai}\,x_{m+c,k}$
this sum equals the sum of the series in the lines 
\eqref{sum1} and \eqref{sum2}.
This equality proves that the map $\chi$ intertwines the actions of
the Yangian $\Y(\gl_n)\ts$.
\qed

\medskip\smallskip
By the transitivity of induction, Theorem \ref{parind}
can be extended from the maximal to all
parabolic subalgebras of the Lie algebra $\gl_m\ts$.
Consider the Borel
subalgebra $\h\op\np$ of $\gl_m\ts$.
Here $\h$ is the Cartan subalgebra of $\gl_m$ spanned
by the elements $E_{aa}\ts$, and
$\n'$ is the nilpotent subalgebra of $\gl_m$ spanned
by the elements $E_{ab}$ with $a<b\ts$.

Take any element $\mu$ of the vector space $\h^\ast$ dual to
$\h\ts$, any such element is called a \textit{weight}\ts.
The weight $\mu$ can be identified with
the sequence $(\ts\mu_1\lcd\mu_m)$ of its \textit{labels},
where $\mu_a=\mu(E_{aa})$ for each $a=1\lcd m\ts$.
Consider the Verma module $M_{\ts\mu}$ over the Lie algebra
$\gl_m\ts$. It can be described as the quotient of the algebra
$\U(\gl_m)$ by the left ideal generated by all the elements
$E_{ab}$ with $a<b$ and the elements $E_{aa}-\mu_a\ts$.
The elements of the Lie algebra $\gl_m$ act on this quotient
via left multiplication. The image of the element $1\in\U(\gl_m)$
in this quotient is denoted by $1_\mu\ts$. Then
$X\cdot 1_\mu=0$ for all $X\in\n'$ while
$$
X\cdot 1_\mu=\mu\ts(X)\cdot 1_\mu
\quad\textrm{for all}\quad
X\in\h\ts.
$$

Let us now
apply the functor \eqref{zelefun} to the $\gl_m$-module $V=M_{\ts\mu}\ts$,
and the functor \eqref{drinfun} to the resulting $\H_N$-module
$$
W=(\ts M_{\ts\mu}\ot(\CC^{\ts m})^{\ot N})_{\ts\n}\ts.
$$
We obtain the $\Y(\gl_n)$-module
$$
((\ts M_{\ts\mu}\ot(\CC^{\ts m})^{\ot N})_{\ts\n}\ot
(\CC^{\ts n})^{\ot N})^{\Sym_N}_{\,-}
=
(\ts M_{\ts\mu}\ot\S^N(\CC^{\ts m}\ot\CC^{\ts n}))_{\ts\n}\ts.
$$
By taking the direct sum over $N=0,1,2,\ts\ldots$
of these $\Y(\gl_n)$-modules, we get the $\Y(\gl_n)$-module
$$
\A_{\ts m}(\ts M_{\ts\mu})_{\ts\n}=
(\ts M_{\ts\mu}\ot\S\ts(\CC^{\ts m}\ot\CC^{\ts n}))_{\ts\n}\,.
$$
Note that $\A_{\ts m}(\ts M_{\ts\mu})_{\ts\n}$ is also a module over
the Cartan subalgebra $\h\ts$. Using the basis
$E_{11}\lcd E_{mm}$ identify $\h$ with the direct sum of $m$
copies of the Lie algebra $\gl_{\ts1}\ts$.
Consider the Verma modules $M_{\ts\mu_1}\lcd M_{\ts\mu_m}$ over
$\gl_{\ts1}\ts$. By applying
Theorem~\ref{parind} we get the next result,
which can be also derived from \cite[Theorem 3.3.1]{AS}.

\begin{corollary}
\label{bimequiv}
The bimodule $\A_{\ts m}(\ts M_{\ts\mu})_{\ts\n}$ of\/ $\h$ and\/ $\Y(\gl_n)$
is equivalent to the tensor product
$$
\A_{\ts1}\ts(\ts M_{\ts\mu_1})
\ot
\A_{\,1}^{\tts1}\ts(\ts M_{\ts\mu_2})
\ot
\ldots
\ot
\A_{\,1}^{\,m-1}\ts(\ts M_{\ts\mu_m})\ts.
$$ 
\end{corollary}

We complete this section with describing 
for any $t\com z\in\CC$ the bimodule $\A_{\ts1}^{\ts z}\ts(M_{\ts t})$ 
over $\gl_{\ts1}$ and $\Y(\gl_n)\ts$.
The Verma module $M_{\ts t}$ over
$\gl_{\ts1}$ is one-dimensional, and the element
$E_{11}\in\gl_{\ts1}$ acts on $M_{\ts t}$ as multiplication by $t\ts$.
The vector space of bimodule 
$\A_{\ts1}(M_{\ts t})$ is the exterior algebra 
$\S\ts(\CC^{\ts1}\ot\CC^{\ts n})=\S\ts(\CC^{\ts n})\ts$,
which we identify with 
$\P\ts(\CC^{\ts1}\ot\CC^{\ts n})=\P\ts(\CC^{\ts n})\ts$. 
Then $E_{11}$ acts on $\A_{\ts1}(M_{\ts t})$ as the differential operator
$$
t\,+\,\sum_{k=1}^n\,x_{1k}\,\d_{\ts1k}\,.
$$
The action of $E_{11}$ on 
$\A_{\ts1}^{\ts z}\ts(M_{\ts t})$ is the same as on $\A_{\ts1}(M_{\ts t})\ts$.
The generator $T_{ij}^{\ts(s+1)}$ of $\Y(\gl_n)$ with $s=0,1,2,\ts\ldots$
acts on $\A_{\ts1}(M_{\ts t})$ as the differential operator 
$
t^{\ts s}\,x_{1i}\,\d_{\ts1j}\,;
$
this is what Proposition \ref{dast} states in the case $m=1$.
Note that the operator $x_{1i}\,\d_{\ts1j}$ describes the
action on $\P\ts(\CC^{\ts1}\ot\CC^{\ts n})=\P\ts(\CC^{\ts n})$ 
of the element $E_{ij}\in\gl_n\ts$.
Hence the action of the algebra $\Y(\gl_n)$ on 
$\A_{\ts1}(M_{\ts w})$ can be obtained from the action
of $\gl_n$ on $\P\ts(\CC^{\ts n})\ts$ by pulling back through 
the homomorphism $\pi_n:\Y(\gl_n)\to\U(\gl_n)\ts$, and then
through the automorphism $\tau_{w}$ of $\Y(\gl_n)\ts$; see
the definitions \eqref{tauz},\eqref{pin}.
Hence the action of $\Y(\gl_n)$ on $\A_{\ts1}^{\ts z}\ts(M_{\ts t})$
can be obtained from the action
of $\gl_n$ on $\P\ts(\CC^{\ts n})\ts$ by pulling back through 
$\pi_n$ and then through the automorphism $\tau_{\ts t-z}\ts$.


\section*{3.\ Zhelobenko operators}
\setcounter{section}{3}
\setcounter{equation}{0}
\setcounter{theorem}{0}

Consider the group $\Sym_m$ as the Weyl group of
the reductive Lie algebra $\gl_m\ts$. This group acts on the 
vector space $\gl_m$ so that for any $\si\in\Sym_m$ and $a,b=1\lcd m$
$$
\si:E_{ab}\mapsto E_{\si(a)\ts\si(b)}\ts.
$$
This action extends to an action of $\Sym_m$
by automorphisms of the associative algebra $\U(\gl_m)\ts$.
The group $\Sym_m$ also acts on the vector space $\h^\ast\ts$.
Let $E_{11}^{\,\ast}\lcd E_{mm}^{\,\ast}$ be the basis
of $\h^\ast$ dual to the basis $E_{11}\lcd E_{mm}$ of $\h\ts$. Then
$$
\si:E_{aa}^{\,\ast}\mapsto E_{\si(a)\ts\si(a)}^{\,\ast}\ts.
$$
If we identify each weight $\mu\in\h^\ast$ with
the sequence $(\ts\mu_1\lcd\mu_m)$ of its labels, then
$$
\si:(\ts\mu_1\lcd\mu_m)\mapsto
(\,\mu_{\ts\si^{-1}(1)}\lcd\mu_{\ts\si^{-1}(m)}\ts)\ts.
$$
Let $\rho\in\h^\ast$ be the weight with
that sequence of labels $(\ts0,-1\lcd1-m)\ts$.
The \textit{shifted\/} action of any element $\si\in\Sym_m$ on $\h^\ast$
is defined by the assignment
\begin{equation}
\label{saction}
\mu\,\mapsto\,\si\circ\mu=\si\ts(\mu+\rho)-\rho\ts.
\end{equation}

For $a,b=1\lcd m$ put $\ep_{ab}=E_{aa}^{\,\ast}-E_{bb}^{\,\ast}\,$.
The elements $\ep_{ab}\in\h^\ast$ with $a<b$ and $a>b$ are
the \textit{positive\/} and \textit{negative roots\/} 
respectively; $\ep_{ab}=0$ when $a=b\ts$.
The elements $\ep_c=\ep_{c,c+1}\in\h^\ast$ with $c=1\lcd m-1$ are the
\textit{simple\/} positive roots. Put
$$
E_{c}=E_{c,c+1}\ts,
\quad
F_{c}=E_{c+1,c}
\quad\text{and}\quad
H_c=E_{cc}-E_{c+1,c+1}\ts.
$$
For any $a=1\lcd m-1$ these three elements of the
Lie algebra $\gl_m$ span a subalgebra 
isomorphic to the Lie algebra $\mathfrak{sl}_{\ts2}\ts$.

For any $\gl_m$-module $V$ and any $\la\in\h^\ast$
a vector $v\in V$ is said to be \textit{of weight} $\la$ if
$X\,v=\la\ts(X)\,v$ for any $X\in\h\ts$.
Denote by $V^\la$ the subspace in $V$ formed by all
vectors of weight $\la\ts$. Recall that $\n$ is
the nilpotent subalgebra of $\gl_m$ spanned
by the elements $E_{ab}$ with $a>b\ts$.  
In this section, we will employ
the general notion of a \textrm{Mickelsson algebra} introduced in \cite{M1}
and developed by D.\ts Zhelobenko \cite{Z}. 
Namely, we will show how this notion gives rise to
a distinguished $\Y(\gl_n)$-intertwining operator
\begin{equation}
\label{distoper}
\A_{\ts m}(\ts M_{\ts\mu})_{\ts\n}^{\ts\la}
\ \to\,
\A_{\ts m}(\ts M_{\ts\si\,\circ\ts\mu})_{\ts\n}^{\,\si\,\circ\ts\la}
\end{equation}
for any element $\si\in\Sym_m$ and any
weight $\mu\in\h^\ast$ such that 
\begin{equation}
\label{muass}
\mu_a-\mu_b\notin\ZZ
\quad\text{whenever}\quad
a\neq b\ts.
\end{equation}
The source and the target vector
spaces in \eqref{distoper} are non-zero only if
all labels of the weight $\la-\mu$ are non-negative integers. Then 
$\la_a-\la_b\notin\ZZ$ whenever $a\neq b\ts$.

We have a representation
$\ga:\,\U(\gl_m)\to\PD\ts(\CC^{\ts m}\ot\CC^{\ts n})$
such that the image $\ga(E_{ab})$ is the operator \eqref{glmact}.
Note that the group $\Sym_m$ acts by automorphisms of the algebra 
$\PD\ts(\CC^{\ts m}\ot\CC^{\ts n})\ts$, so that for $k=1\lcd n$
\begin{equation}
\label{symact}
\si:\
x_{ak}\ts\mapsto\ts x_{\si(a)\ts k}\,,\
\d_{\ts bk}\ts\mapsto\d_{\ts\si(b)\ts k}\,.
\end{equation}
The homomorphism $\ga$ is $\Sym_m$-equivariant. Denote by $\Ar$
the associative algebra generated by the algebras
$\U(\gl_m)$ and $\PD\ts(\CC^{\ts m}\ot\CC^{\ts n})$
with the cross relations
\begin{equation}
\label{defar}
[X\com Y]=[\ts\ga(X)\com Y]
\end{equation}
for any $X\in\gl_m$ and $Y\in\PD\ts(\CC^{\ts m}\ot\CC^{\ts n})\ts$.
Here each pair of the square brackets denotes the commutator in $\Ar\ts$.
The algebra $\Ar$ is isomorphic to the 
the tensor product of associative algebras \eqref{tenprod}.
The isomorphism can be defined by mapping
the elements $X\in\gl_m$ and $Y\in\PD\ts(\CC^{\ts m}\ot\CC^{\ts n})$
of $\Ar$ respectively to the elements 
$$
X\ot1+1\ot\ga(X)
\quad\textrm{and}\quad
1\ot Y
$$
of \eqref{tenprod}.
This isomorphism is $\Sym_m$-equivariant, and the image 
of the element $E_{ab}\in\gl_m$ under this isomorphism
equals \eqref{eabact}. We will use this isomorphism later on.

Denote by $\J$ the right ideal of the algebra $\Ar$
generated by the elements of the subalgebra $\n\subset\gl_m\ts$.
Let ${\rm Norm}\ts(\J)\ts\subset\ts\Ar$
be the normalizer of this right ideal, so that
$Y\in{\rm Norm}\ts(\J)$ if and only if $\,Y\ts\J\subset\J\ts$.
Then $\J$ is a two-sided ideal of ${\rm Norm}\ts(\J)\ts$.
Our particular \textit{Mickelsson algebra\/} is the quotient 
\begin{equation}
\label{malg}
\R\,=\,\J\,\backslash\,{\rm Norm}\ts(\J)\ts.
\end{equation}


\medskip
\noindent\textit{Remark.}
Via its isomorphism with \eqref{tenprod},
the associative algebra $\Ar$ acts on the tensor product
$V\ot\ts\P\ts(\CC^{\ts m}\ot\CC^{\ts n})$
for any $\gl_m$-module $V$.
The defining embedding of $\gl_m$ into $\Ar$ corresponds to the diagonal action
of the Lie algebra $\gl_m$ on this tensor product.
The Mickelsson algebra $\R$ then acts
on the space of $\n$-coinvariants of $\gl_m$-module 
$V\ot\ts\P\ts(\CC^{\ts m}\ot\CC^{\ts n})\ts$.
\qed

\medskip\smallskip
Let $\,\overline{\!\U(\h)\!\!\!}\,\,\,$ 
be the ring of fractions of 
$\U(\h)$ relative to the set of denominators
\begin{equation}
\label{denset}
\{\,E_{aa}-E_{bb}+z\ |\ 1\le a\com b\le m\ts;\ a\neq b\ts;\ z\in\ZZ\,\ts\}\,.
\end{equation}
The elements of this ring can also
be regarded as rational functions on the vector space
$\h^\ast\ts$. The elements of 
$\U(\h)\subset\,\overline{\!\U(\h)\!\!\!}\,\,\,$
are then regarded as polynomial functions. 
Denote by $\Ab$ the ring of fractions of $\Ar$ 
relative to the same set of denominators \eqref{denset},
regarded as elements of $\Ar$ using the embedding of $\h\subset\gl_m$ 
into $\Ar\ts$. The ring $\Ab$ is defined 
due to relations
in $\U(\gl_m)$ and $\Ar\ts$: for $a,b=1\lcd m$
and any $H\in\h$
$$
[\ts H\ts,E_{ab}\ts]=\ep_{ab}(H)\ts E_{ab}\,,
\quad
[\ts H\,,x_{ak}\ts]=
E_{aa}^{\,\ast}(H)\,x_{ak}\,,
\quad
[\ts H\,,\d_{\ts bk}\ts]=-\ts
E_{bb}^{\,\ast}(H)\,\d_{\ts bk}\,.
$$
Therefore the ring $\Ar$ satisfies the Ore condition relative
to its subset \eqref{denset}.
Using left multiplication by elements of
$\,\overline{\!\U(\h)\!\!\!}\,\,\,$,
the ring $\Ab$ becomes a module over
$\,\overline{\!\U(\h)\!\!\!}\,\,\,\ts$.

The ring $\Ab$ is also an associative algebra over the field $\CC\ts$.
For $c=1\lcd m-1$ define a linear map $\xi_{\ts c}:\Ar\to\Ab$
by setting
\begin{equation}
\label{q1}
\xi_{\ts c}\,(\ts Y)=Y+\,\sum_{s=1}^\infty\,\,
(\ts s\ts!\,H_c^{\ts(s)}\ts)^{-1}\ts E_c^{\ts s}\,
\widehat{F}_c^{\ts s}(\ts Y)
\end{equation}
for $Y\in\Ar\ts$. Here
$$
H_c^{\ts(s)}=H_c(H_c-1)\ldots(H_c-s+1)
$$
and $\widehat{F}_c$ is the operator of adjoint action
corresponding to the element $F_c\in\Ar\ts$,
$$
\widehat{F}_c(\ts Y)=[\ts F_c\ts\com Y\ts]\ts.
$$
For any given element $Y\in\Ar$ only finitely many terms of the sum 
\eqref{q1} differ from zero, hence the map $\xi_{\ts c}$ is well defined.
The definition \eqref{q1} and the following proposition
go back to \cite[Section 2\ts]{Z}.
Put $\Jb=\,\overline{\!\U(\h)\!\!\!}\,\,\,\,\J\ts$.
Then $\Jb$ is a right ideal of $\Ab\ts$.

\begin{proposition}
\label{prop1}
For any $X\in\h$ and $Y\in\Ar$ we have
\begin{equation}
\label{q11}
\xi_{\ts c}(X\ts Y)\in
(\ts X+\ep_c(X))\,\ts\xi_{\ts c}(\ts Y)\ts+\ts\Jb\ts,
\end{equation}
\begin{equation}
\label{q12}
\xi_{\ts c}(\ts Y X)\in\,
\xi_{\ts c}(\ts Y)\ts(\ts X+\ep_c(X))\ts+\ts\Jb\ts.
\end{equation}
\end{proposition}

This proposition coincides with \cite[Proposition 3.1]{KN1} 
so we skip the proof here.
The property \eqref{q11} allows us to define a linear map 
$
\bar\xi_{\ts c}:\Ab\to\Jb\,\ts\backslash\,\Ab
$
by setting
$$
\bar\xi_{\ts c}(X\,Y)=Z\,\xi_{\ts c}(\ts Y)\ts+\ts\Jb
\quad\ \text{for any}\quad
X\in\,\overline{\!\U(\h)\!\!\!}\,\,\,
\quad\text{and}\quad
Y\in\Ar\ts,
$$
where the element $Z\in\,\overline{\!\U(\h)\!\!\!}\,\,\,$
is defined by the equality
$$
Z(\mu)=X(\ts\mu+\ep_c)
\quad\ \text{for any}\quad
\mu\in\h^\ast
$$
when $X$ and $Z$ are
regarded as rational functions on $\h^\ast\ts$.

The action of the group $\Sym_m$ on the algebra $\U(\gl_m)$ 
extends to an action on $\,\overline{\!\U(\h)\!\!\!}\,\,\,\ts$, so that
for any $\si\in\Sym_m$
$$
(\ts\si\ts X)(\mu)=X(\ts\si^{-1}(\mu))
$$
when the element $X\in\,\overline{\!\U(\h)\!\!\!}\,\,\,$
is regarded as a rational function on $\h^\ast\ts$. The action
of $\Sym_m$ by automorphisms of the algebra $\Ar$ then 
extends to an action by automorphisms of $\Ab\ts$.
For any $c=1\lcd m-1$ let $\si_c\in\Sym_m$ be
the transposition of $c$ and $c+1\ts$.
Consider the image $\si_c(\ts\Jb\ts)\ts$,
this is again a right ideal of $\Ab$. Next proposition
also goes back to \cite{Z}. 
It coincides with \cite[Proposition 3.2]{KN1} 
so we skip the proof~here.

\begin{proposition}
\label{prop2}
We have\/ $\si_c(\ts\Jb\ts)\subset\ker\ts\bar\xi_{\ts c}\ts$.
\end{proposition}

This proposition allows us to define for any $c=1\lcd m-1$ a linear map
\begin{equation}
\label{xic}
\xic_{\ts c}:\,\Jb\,\ts\backslash\,\Ab\to\Jb\,\ts\backslash\,\Ab
\end{equation}
as the composition $\bar\xi_{\ts c}\,\si_c$ 
applied to the elements of $\Ab$ which were
taken modulo $\Jb\ts$. 

\medskip\smallskip
\noindent\textit{Remark.}
Observe that $\U(\h)\ts\subset\ts{\rm Norm}\ts(\J)\ts$.
Let us denote by $\,\overline{\!{\rm Norm}\ts(\J)\!\!\!}\,\,\,$
the ring of fractions of ${\rm Norm}\ts(\J)$
relative to the same set of denominators \eqref{denset} as before.
Evidently, then $\Jb$ is a two-sided ideal of the ring
$\,\overline{\!{\rm Norm}\ts(\J)\!\!\!}\,\,\,\ts$.
The quotient ring
$$
\Rb\,=\,\Jb\ts\,\backslash\,\,\overline{\!{\rm Norm}\ts(\J)\!\!\!}\,\,\,
$$
bears the same name of \textrm{Mickelsson algebra\/},
as the quotient ring \eqref{malg} does.
One can show \cite{KO} that the linear map \eqref{xic} preserves
the subspace $\Rb\subset\Jb\,\ts\backslash\,\Ab\ts$,
and determines an automorphism of the algebra $\Rb\ts$.
We do not use these two facts, but
our construction of the $\Y(\gl_n)$-intertwining operator
\eqref{distoper} is underlied by them.
\qed

\medskip\smallskip
In their present form, the operators $\xic_{\ts1}\lcd\xic_{\ts m-1}$  
on the vector space $\Jb\,\ts\backslash\,\Ab$
were introduced in \cite{KO}.
We will call them \textit{Zhelobenko operators}.
The next proposition states the key property of these
operators; for the proof see \cite[Section 6]{Z}.

\begin{proposition}
The operators\/ $\xic_{\ts1}\lcd\xic_{\ts m-1}$ 
on\/ $\Jb\,\ts\backslash\,\Ab$ satisfy the braid relations
$$
\begin{aligned}
\xic_{\ts c}\,\xic_{\ts c+1}\,\xic_{\ts c}
&\,=\,
\xic_{\ts c+1}\,\xic_{\ts c}\,\xic_{\ts c+1}
\,\quad\quad\textit{for}\ \quad
c<m-1\ts,
\\
\xic_{\ts b}\,\xic_{\ts c}
&\,=\,
\xic_{\ts c}\,\xic_{\ts b}
\hspace{39pt}
\ \quad\textit{for}\ \quad
b<c-1\ts.
\end{aligned}
$$
\end{proposition}

\begin{corollary}
\label{decindep}
For any reduced decomposition $\si=\si_{c_1}\ldots\ts\si_{c_K}$ 
in the group\/ $\Sym_m$ the composition
$\xic_{\ts c_1}\ldots\,\xic_{\ts c_K}$
of operators on\/ $\Jb\,\ts\backslash\,\Ab$ 
does not depend on the choice of the decomposition of\/ $\si\ts$.
\end{corollary}

Recall that $\n'$ denotes the nilpotent subalgebra of $\gl_m$ spanned
by the elements $E_{ab}$ with $a<b\ts$.
Denote by $\Jp$ the left ideal of the algebra $\Ar$ generated by the
elements of the subalgebra $\np\subset\gl_m\ts$. 
Put $\Jpb=\,\overline{\!\U(\h)\!\!\!}\,\,\,\,\Jp\ts$.
Then $\Jpb$ is a left ideal of $\Ab\ts$.
Consider the image $\si_c(\ts\Jpb\ts)\ts$,
this is again a left ideal of $\Ab$.

\begin{proposition}
\label{prop3}
We have\/ $\bar\xi_{\ts c}(\,\si_c(\Jpb\ts))\ts\subset\ts\Jpb\ns+\ts\Jb\ts$.
\end{proposition}

This proposition coincides with \cite[Proposition 3.5]{KN1} 
so we again skip the proof. 
Proposition \ref{prop3} implies that for each $c=1\lcd m-1$
the Zhelobenko operator \eqref{xic} induces a linear map
$$
\Jb\,\ts\backslash\,\Ab\,/\,\Jpb\,\to\,\Jb\,\ts\backslash\,\Ab\,/\,\Jpb\,.
$$
Now take a weight $\mu\in\h^\ast$ satisfying \eqref{muass}. 
We shall keep the assumption \eqref{muass} on $\mu$
till the end of this section.
Let $\I_{\ts\mu}$ be the left ideal of the algebra $\Ar$
generated by the elements 
$$
E_{ab}
\quad\text{with}\quad
a<b\ts,
\quad
E_{aa}-\mu_a\ts
\,\quad\text{and}\,\quad
\d_{\ts bk}
$$
for all possible $a,b$ and $k\ts$.
Under the isomorphism of $\Ar$ 
with the tensor product \eqref{tenprod}, 
the ideal $\I_{\ts\mu}$ of $\Ar$ corresponds to the ideal of 
\eqref{tenprod} generated by the elements 
$$
E_{ab}\ot1
\quad\text{with}\quad
a<b\ts,
\quad
E_{aa}\ot1-\mu_a\ts
\,\quad\text{and}\,\quad
1\ot\d_{\ts bk}
$$
for all possible $a,b$ and $k\ts$.
Indeed, for any $a,b=1\lcd m$
the image of the element $E_{ab}\in\Ar$ in the algebra \eqref{tenprod}
is the sum \eqref{eabact}, which equals $E_{ab}\ot1$ plus elements
divisible on the right by the tensor products of the form
$1\ot\d_{\ts bk}\ts$. But the quotient space of \eqref{tenprod}
with respect to the latter ideal can be naturally identified
with the tensor product
$M_{\ts\mu}\ot\P(\ts\CC^{\ts m}\ot\CC^{\ts n})\ts$.
Using the isomorphism of algebras $\Ar$ and \eqref{tenprod},
the quotient space $\Ar\,/\,\I_{\ts\mu}$ can be then also identified with
$M_{\ts\mu}\ot\P(\ts\CC^{\ts m}\ot\CC^{\ts n})\ts$.

Note that
$\mu(H_c)\notin\ZZ$ for any index $c=1\lcd m-1$ due to \eqref{muass}.
Hence we can define the subspace
\,$\Ib_{\ts\mu}=\,\overline{\!\U(\h)\!\!\!}\,\,\,\,\I_{\ts\mu}$ of $\Ab\ts$.
This subspace is also a left ideal of the algebra $\Ab\ts$.
The quotient space $\Ab\,/\,\Ib_{\ts\mu}$ can be still identified 
with the tensor product
$M_{\ts\mu}\ot\P(\ts\CC^{\ts m}\ot\CC^{\ts n})\ts$.
The quotient of $\Ab$ by $\Ib_{\ts\mu}$ and $\Jb$ 
can be then identified with the space of $\n$-coinvariants,
\begin{equation}
\label{dqci}
\Jb\,\ts\backslash\,\Ab\,/\,\Ib_{\ts\mu}
\,=
(\ts M_{\ts\mu}\ot\P(\ts\CC^{\ts m}\ot\CC^{\ts n}))_{\ts\n}\ts.
\end{equation}

Consider the left ideal of the algebra $\Ar$ generated by all the elements
$\d_{\ts bk}$ where $b=1\lcd m$ and $k=1\lcd n\ts$. 
By the definition \eqref{q1}, the image of this ideal under the
map $\xi_{\ts c}$ is contained in the left ideal of $\Ab$
generated by the same elements. The latter ideal is preserved
by the action on $\Ab$ of the element $\si_c\in\Sym_m\ts$. By \eqref{saction},
$$
\si_c(\ts\mu+\ep_c)=\si_c\circ\mu\ts.
$$
The property \eqref{q12} and Proposition \ref{prop3} imply that 
the Zhelobenko operator \eqref{xic} induces a linear map
$$
\Jb\,\ts\backslash\,\Ab\,/\,\Ib_{\ts\mu}
\,\to\,
\Jb\,\ts\backslash\,\Ab\,/\,\Ib_{\,\ts\si_c\ts\circ\ts\mu}\,.
\vspace{2pt}
$$
Via the identifications 
\eqref{dqci}, 
the Zhelobenko operator \eqref{xic} induces a linear map
\begin{equation}
\label{indmap}
(\ts M_{\ts\mu}\ot\P(\ts\CC^{\ts m}\ot\CC^{\ts n}))_{\ts\n}
\,\to\,
(\ts M_{\ts\si_c\ts\circ\ts\mu}\ot\P(\ts\CC^{\ts m}\ot\CC^{\ts n}))_{\ts\n}
\,.
\vspace{4pt}
\end{equation}

\begin{proposition}
\label{inserted}
For any $s=0,1,2,\ts\ldots$
the map \eqref{indmap} commutes with the action of generator
$T_{ij}^{(s+1)}$ of\/ $\Y(\gl_n)$ on the source
and target vector spaces as \eqref{combact}. 
\end{proposition}

\begin{proof}
Let $Y$ be the element of the algebra $\Ar$
corresponding to the element $\eqref{combact}$ of
the algebra \eqref{tenprod} under the isomorphism of these
two algebras. The element $Y$ then belongs to the 
centralizer of the subalgebra $\U(\gl_n)$ in $\Ar\ts$.
So the left multiplication in $\Ab$ by $Y$
preserves the right ideal \hbox{\ts$\Jb\subset\Ab\ts$,}
and commutes with the linear map 
$\bar\xi_{\ts c}:\Ab\to\Jb\,\ts\backslash\,\Ab\,$; 
see the definition \eqref{q1}.
This left multiplication also commutes with the action of the
element $\si_c\in\Sym_m$ on $\Ab\ts$, because 
the element $Y$ of the algebra $\Ar$ is $\Sym_m$-invariant.
\end{proof}

The property \eqref{q11} implies that
the restriction of the linear map \eqref{indmap} to the subspace
of vectors of weight $\la$ is a map
$$
(\ts M_{\ts\mu}\ot\P(\ts\CC^{\ts m}\ot\CC^{\ts n}))_{\ts\n}^{\ts\la}
\,\,\to\,
(\ts M_{\ts\si_c\ts\circ\ts\mu}\ot\P(\ts\CC^{\ts m}
\ot\CC^{\ts n}))_{\ts\n}^{\,\si_c\,\circ\ts\la}\,.
$$
Denote the latter by $I_{\ts c}\,$,
it commutes with the action of $\Y(\gl_n)\ts$
by Proposition \ref{inserted}.

By choosing a reduced decomposition $\si=\si_{c_1}\ldots\ts\si_{c_K}$ 
and taking the composition of operators $I_{c_1}\ldots\ts I_{c_K}$
we obtain an $\Y(\gl_n)$-intertwining operator
$$
I_{\ts\si}:\,
(\ts M_{\ts\mu}\ot\P\ts(\ts\CC^{\ts m}\ot\CC^{\ts n}))_{\ts\n}^{\ts\la}
\,\,\to\,
(\ts M_{\ts\si\,\circ\ts\mu}\ot\P\ts(\ts\CC^{\ts m}
\ot\CC^{\ts n}))_{\ts\n}^{\,\si\,\circ\ts\la}\,.
$$
It does not depend on the choice of the decomposition of
$\si\in\Sym_m$ by Corollary~\ref{decindep}.
This is the operator \eqref{distoper} 
which we intended to exhibit.
Here we identified the exterior algebra 
$\S\ts(\ts\CC^{\ts m}\ot\CC^{\ts n})$ with the ring
$\P\ts(\ts\CC^{\ts m}\ot\CC^{\ts n})$ 
as we did in Section~1.

{}From now on we will assume that all labels of the weight $\nu=\la-\mu$
belong to the set $\{\ts0\com1\lcd n\ts\}\,$;
otherwise both the source and
target modules in \eqref{distoper} are zero. Let
$(\ts\nu_1\lcd\nu_m)$ be the sequence of these labels. 
Consider the element
$$
w_\nu=
(\ts x_{11}\ldots\,x_{1\nu_1}\ts)
\,\ldots\,
(\ts x_{m1}\ldots\,x_{m\nu_m}\ts)
$$
of $\P\ts(\ts\CC^{\ts m}\ot\CC^{\ts n})\ts$.
Note that $w_\nu$
is a \textit{highest\/} vector with respect to the 
action of $\gl_n$ on
$\P\ts(\ts\CC^{\ts m}\ot\CC^{\ts n})\ts$,
any element $E_{ij}\in\gl_n$ with $i<j$
acts on this vector as zero. With respect to the 
action of $\gl_m\ts$,
the vector $w_\nu$ is of weight $\nu\ts$.
Then consider
\begin{equation}
\label{munuvec}
1_\mu\ot w_\nu
\,\in\ts
M_{\ts\mu}\ot\P\ts(\ts\CC^{\ts m}\ot\CC^{\ts n})
\end{equation}
With respect to the action of the Lie algebra $\gl_m$
the vector \eqref{munuvec} is of weight $\la\ts$. 
Denote by $v_{\ts\mu}^{\ts\la}$ the image of the vector \eqref{munuvec} in 
$(\ts M_{\ts\mu}\ot\P(\ts\CC^{\ts m}\ot\CC^{\ts n}))_{\ts\n}^{\ts\la}\,$.

\begin{proposition}
\label{isis}
Under condition \eqref{muass}, the vector
$I_{\ts\si}(\ts v_{\ts\mu}^{\ts\la}\ts)$ equals\/
$v_{\ts\si\,\circ\ts\mu}^{\ts\si\,\circ\ts\la}$ times
\begin{equation}
\label{isim}
\prod\limits_{\substack{1\le a<b\le m \\ \si(a)>\si(b)}}
(-1)^{\ts\nu_a\nu_b}\ 
\left\{
\begin{array}{ll}
\displaystyle
\frac{\,\la_a-\la_b-a+b}{\,\mu_a-\mu_b-a+b}
&\ \textrm{if}\ \ \nu_a<\nu_b\ts;\\[12pt]
\hspace{33pt}1
&\ \textrm{if}\ \ \nu_a\ge\nu_b\ts.
\end{array}
\right.   
\end{equation}

\end{proposition}

\begin{proof}
It suffices to prove this for $\si=\si_c$ with $c=1\lcd m-1\ts$.
Moreover, it then suffices to consider only the case when $m=2$ and hence
$c=1\ts$. Suppose this is the case.
Then under the identification \eqref{dqci} the vector
$$
v_{\ts\mu}^{\ts\la}
\,\in\ts
(\ts M_{\ts\mu}\ot\P(\ts\CC^{\ts 2}\ot\CC^{\ts n}))_{\ts\n}
$$
gets identified with the image in the quotient space
$\Jb\,\ts\backslash\,\Ab\,/\,\Ib_{\ts\mu}$ of the element
$$
w_\nu=(\ts x_{11}\ldots\,x_{1\nu_1}\ts)\,(\ts x_{21}\ldots\,x_{2\nu_2}\ts)
\,\in\ts\Ar\subset\Ab\ts.
$$
According to \eqref{symact},
by applying the transposition $\si_1\in\Sym_2$ to $w_\nu$
we obtain
$$
w=(\ts x_{21}\ldots\,x_{2\nu_1}\ts)\,(\ts x_{11}\ldots\,x_{1\nu_2}\ts)
\,\in\ts\Ar\,.
$$
Note that the vector
$$
v_{\ts\si_1\,\circ\ts\mu}^{\ts\si_1\,\circ\ts\la}
\,\in\ts
(\ts M_{\ts\si_1\ts\circ\ts\mu}\ot\P(\ts\CC^{\ts 2}\ot\CC^{\ts n}))_{\ts\n}
$$
is identified with the image in the quotient space
$
\Jb\,\ts\backslash\,\Ab\,/\,\Ib_{\,\ts\si_1\ts\circ\ts\mu}
$
of the element of $\Ar\ts$,
$$
w_{\ts\si_1(\nu)}=
(\ts x_{11}\ldots\,x_{1\nu_2}\ts)\,(\ts x_{21}\ldots\,x_{2\nu_1}\ts)
=(-1)^{\ts\nu_1\nu_2}\,w\ts.
$$

By applying the map $\xi_{1}$ to the element $w\in\Ar$ we get
the sum of elements of $\Ab\ts$,
\begin{equation}
\label{sumel}
\sum_{s=0}^\infty\,\,
(\ts s\ts!\,H_1^{\ts(s)}\ts)^{-1}\ts E_1^{\ts s}\,
\widehat{F}_1^{\ts s}\ts(\ts w\ts)\,.
\end{equation}
In particular, here by the definition of the algebra $\Ar$ we have
$$
\widehat{F}_1(\ts w\ts)=
[\,F_1\com\ts w\,]=
\sum_{k=1}^n\,\ts
[\,x_{2k}\,\d_{\ts1k}\ts\com\ts w\,]\,.
$$
If $\nu_1\ge\nu_2$
then ${F}_1(\ts w\ts)=0\,$ and $\,\xic_1(\ts w_\nu\ts)=\xi_1(\ts w\ts)=w\ts$,
as required in this case.

Now suppose that $\nu_1<\nu_2\ts$. Denote $d=\nu_2-\nu_1\ts$. Then 
$\widehat{F}_1^{\ts s}\ts(\ts w\ts)$ equals the sum over all subsets 
$\mathcal{I}_{\ts s}\subset\{\ts\nu_1+1\lcd\nu_2\ts\}$ with cardinality 
$s\ts$, of the elements
$$
\wp=s\ts!\,
(\ts x_{21}\ldots\,x_{2\nu_1}\ts)\,
(\ts x_{11}\ldots\,x_{1\nu_1}\ts)\,
(\ts x_{c_1\nu_1+1}\ldots\,x_{c_d\nu_2}\ts)
$$
where for $k=1\lcd d$ we have $c_k=2$ or $c_k=1$
depending on whether $\nu_1+k\in\mathcal{I}_{\ts s}$ or not.
Since the element $E_1\in\Ar$
is a generator of the left ideal $\I_{\,\ts\si_1\ts\circ\ts\mu}\ts$,
for any element $Y\in\PD\ts(\ts\CC^{\ts 2}\ot\CC^{\ts n})$ 
the image of the product $E_1\ts Y$ in the quotient vector space
$\Jb\,\ts\backslash\,\Ab\,/\,\Ib_{\,\ts\si_1\ts\circ\ts\mu}$
coincides with the image of
$$
[\,E_1\com\ts Y\,]=
\sum_{k=1}^n\,\ts
[\,x_{1k}\,\d_{\ts2k}\ts\com\ts Y\,]\,.
$$
It follows that for the summand $\wp$ corresponding to any subset 
$\mathcal{I}_{\ts s}\ts$, the image of the product $E_1^{\ts s}\,\wp$ in
$\Jb\,\ts\backslash\,\Ab\,/\,\Ib_{\,\ts\si_1\ts\circ\ts\mu}$
coincides with the image of $(\ts s\ts!\ts)^{\ts2}\,w\ts$. 
Since the total number of
the subsets $\mathcal{I}_{\ts s}$ is
$d\,!\,/\ts s\ts!\,(d-s)\ts!\,$ for any given $s\le d\ts$, 
the image of the sum \eqref{sumel} in 
$\Jb\,\ts\backslash\,\Ab\,/\,\Ib_{\,\ts\si_1\ts\circ\ts\mu}$
coincides with that of the sum
\begin{equation}
\label{sumelll}
\sum_{s=0}^d\,\,
d\,!\,(\ts(d-s)\ts!\,H_1^{\ts(s)}\ts)^{-1}\ts w\ts.
\end{equation}

In the sum \eqref{sumelll} the symbol $H_1^{\ts(s)}$ 
stands for the product in the algebra $\Ar\ts$,
$$
\prod_{r=1}^s\,
(\ts H_1-r+1\ts)
\,=\ts
\prod_{r=1}^s\,
(\ts E_{11}-E_{22}-r+1\ts)
\,.
$$
Since the elements
$$
E_{11}-\mu_2+1
\ \quad\textrm{and}\ \quad
E_{22}-\mu_1-1
\vspace{4pt}
$$
are also generators of the left ideal
$\I_{\,\ts\si_1\ts\circ\ts\mu}\ts$,
the image of the product $H_1\ts w$ in the quotient vector space
$\Jb\,\ts\backslash\,\Ab\,/\,\Ib_{\,\ts\si_1\ts\circ\ts\mu}$
coincides with the image of
$$
[\,H_1\com\ts w\,]+(\mu_2-\mu_1-2)\,w\,=
\sum_{k=1}^n\,\ts
[\,x_{1k}\,\d_{\ts1k}-x_{2k}\,\d_{\ts2k}\ts\com\ts w\,]+(\mu_2-\mu_1-2)\,w=
$$
$$
(\nu_2-\nu_1+\mu_2-\mu_1-2)\,w=(\la_2-\la_1-2)\,w\,.
\vspace{6pt}
$$
Therefore
the image of the sum \eqref{sumelll} in
$\Jb\,\ts\backslash\,\Ab\,/\,\Ib_{\,\ts\si_1\ts\circ\ts\mu}$
coincides with the image of
$$
\sum_{s=0}^{d}\,
\,\prod_{r=1}^s\,
\frac{d-r+1}{\ts\la_2-\la_1-r-1\ts}
\,\cdot\,w\,=\, 
\prod_{r=1}^{d}\,
\frac{\ts\la_1-\la_2+d-r+1\ts}{\la_1-\la_2+r+1}
\,\cdot\,w\,=\,
\frac{\ts\la_1-\la_2+1\ts}{\mu_1-\mu_2+1}
\,\cdot\,w
\vspace{2pt}
$$
as required. Here for $t\ts=\la_1-\la_2\ts$ and 
any positive integer $d$ we used the equality of rational functions
of the variable $t\ts$, 
$$
\sum_{s=0}^{d}\,\,
\,\prod_{r=1}^s\,\,
\frac
{d-r+1}{\ts-\,t-r-1\ts}
\,\,\,=\,\,
\frac{\,t+1}{\,t+d+1}
\vspace{2pt}
$$
which can be easily proved by induction on $d=0\com1\com2\com\ts\ldots\,\,$.
\end{proof}

Note that 
the $\Y(\gl_n)$-intertwining operator $I_{\ts\si}$ has been defined only
when the weight $\mu$ satisfies the condition \eqref{muass}. 
We also assume that all labels of the weight
$\nu=\la-\mu$ are non-negative integers.
Recall that the sequence of labels $(\ts\rho_1\lcd\rho_m)$
of the weight $\rho$ is $(\ts0\com-1\lcd1-m\ts)\ts$.
For any $z\in\CC$ and $N=0,1\lcd n$ let us denote by 
$A_{\ts z}^{\ts N}$
the $\Y(\gl_n)$-module obtained from the standard action
of $\U(\gl_n)$ on $\S^N(\CC^{\ts n})\ts$ by pulling back through 
the homomorphism $\pi_n:\Y(\gl_n)\to\U(\gl_n)\ts$, and then
through the automorphism $\tau_z$ of $\Y(\gl_n)\ts$; see
the definitions \eqref{tauz} and \eqref{pin}.
Using Corollary \ref{bimequiv} and the subsequent remarks, 
we can replace the source and target modules in \eqref{distoper}
by equivalent $\Y(\gl_n)$-modules to get an intertwining
operator between two tensor products of $\Y(\gl_n)$-modules,
\begin{equation}
\label{swnop}
A_{\ts\mu_1+\rho_1}^{\ts\nu_1}
\ot\ldots\ot 
A_{\ts\mu_m+\rho_m}^{\ts\nu_m}
\to\,
A_{\ts\widetilde\mu_1+\widetilde\rho_1}^{\ts\widetilde\nu_1}
\ot\ldots\ot 
A_{\ts\widetilde\mu_m+\widetilde\rho_m}^{\ts\widetilde\nu_m}
\end{equation}
where for $a=1\lcd m$ we write
$$
\widetilde\mu_a=\mu_{\si^{-1}(a)}\ts,
\quad 
\widetilde\nu_a=\nu_{\si^{-1}(a)}\ts,
\quad
\widetilde\rho_a=\rho_{\si^{-1}(a)}
$$

It is well known that under the condition \eqref{muass} on
the sequence $(\ts\mu_1\lcd\mu_m)$
both tensor products are irreducible $\Y(\gl_n)$-modules,
equivalent to each other\ts; see
\cite[Theorem 3.4]{NT1}. So an intertwining operator 
between these two tensor products is unique up to a multiplier from
$\CC\ts$. For the operator
corresponding to $I_{\ts\si}$ this multiplier
is determined by Proposition \ref{isis}.
Another expression for an intertwining operator 
between two tensor products of $\Y(\gl_n)$-modules \eqref{swnop}
can be obtained by using a method of I.\,Cherednik \cite{C1},
see for instance \cite[Section 2]{NT1}.

\medskip\smallskip
\noindent\textit{Remark.}
The product \eqref{isim} in Proposition \ref{isis}
does not depend on the choice of reduced decomposition 
$\si_{c_1}\ldots\ts\si_{c_K}$ of $\si\in\Sym_m\ts$.
The uniqueness of the
intertwining operator \eqref{swnop} thus provides another proof
of the independence of the composition
$I_{c_1}\ldots\ts I_{c_K}$ on the choice of the
decomposition of $\si\ts$, not involving Corollary \ref{decindep}.
\qed

\section*{4.\ Olshanski homomorphism}
\setcounter{section}{4}
\setcounter{equation}{0}
\setcounter{theorem}{0}

Let $l$ be a positive integer.
The decomposition $\CC^{\ts n+l}=\CC^{\ts n}\op\CC^{\ts l}$ 
defines an embedding
of the direct sum $\gl_n\op\gl_{\ts l}$ of Lie algebras into $\gl_{n+l}\ts$.
As a subalgebra of $\gl_{n+l}\ts$,
the direct summand $\gl_n$ is spanned by the matrix units 
$E_{ij}\in\gl_{n+l}$ where
$i,j=1\lcd n\ts$. The direct summand $\gl_{\ts l}$ is spanned by
the matrix units $E_{ij}$ where $i,j=n+1\lcd n+l\ts$.
Let $\Cr_{\ts l}$ be the centralizer in 
$\U(\gl_{n+l})$ of the subalgebra $\gl_{\ts l}\subset\gl_{n+l}\ts$.
Set $\Cr_0=\U(\gl_n)\ts$.

Proposition \ref{dast} shows that for any positive integer $m$ a homomorphism
\begin{equation}
\label{lastrem}
\Y(\gl_n)\,\to\ts\U(\gl_m)\ot\PD\ts(\CC^{\ts m}\ot\CC^{\ts n})
\end{equation}
of associative algebras can be defined
by mapping $T_{ij}^{\ts(s+1)}$ to the sum \eqref{combact}.
The image of this homomorphism is contained in the
centralizer of the image of $\gl_m$ in 
$\U(\gl_m)\ot\PD\ts(\CC^{\ts m}\ot\CC^{\ts n})\ts$, 
see the remark just after the proof of Proposition \ref{dast}.
In this section we compare this homomorphism with
a homomorphism $\Y(\gl_n)\to\Cr_{\ts l}$ defined by G.\,Olshanski \cite{O1}.

Consider the Yangian $\Y(\gl_{n+l})\ts$. 
The subalgebra in $\Y(\gl_{n+l})$ generated by 
$$
T_{ij}^{\ts(1)},T_{ij}^{\ts(2)},\ts\ldots
\quad\text{where}\quad
i\ts,\ns j=1\lcd n 
$$
is isomorphic to $\Y(\gl_n)\ts$ as an associative algebra, see
\cite[Corollary 1.23]{MNO}. Thus we have a natural embedding
$\Y(\gl_n)\to\Y(\gl_{n+l})\ts$, which will be denoted by $\io_{\ts l}\ts$.
Note that $\io_{\ts l}$ is not a Hopf algebra homomorphism.
We also have a surjective homomorphism of associative algebras
$$
\pi_{n+l}:\ts\Y(\gl_{n+l})\to\U(\gl_{n+l})\ts,
$$
see \eqref{pin}. The composition 
$\pi_{n+l}\,\io_{\ts l}$ coincides with the homomorphism $\pi_n\ts$.

Further, consider the involutive automorphism $\om_{n+l}$
of the algebra $\Y(\gl_{n+l})\ts$, see the definition \eqref{1.51}.
The image of the composition of homomorphisms
$$
\pi_{n+l}\,\ts\om_{n+l}\,\ts\io_{\ts l}:\ts
\Y(\gl_n)\to
\U(\gl_{n+l})
$$
belongs to the subalgebra $\Cr_{\ts l}\subset\U(\gl_{n+l})\ts$.
This image and the centre of the algebra
of $\U(\gl_{n+l})$ generate the subalgebra $\Cr_{\ts l}\ts$.
For the proofs of these two assertions see \cite[Section 2.1]{O2}.
The \textit{Olshanski homomorphism\/} $\Y(\gl_n)\to\Cr_{\ts l}$
is the composition
\begin{equation}
\label{al}
\al_{\ts l}\ts=\ts
\pi_{n+l}\,\ts
\om_{n+l}\,\ts
\io_{\ts l}\,\ts
\tau_{-\ts l}\ts.
\end{equation}
Set $\al_0=\pi_n\ts$. Further comments on
the family of homomorphisms $\al_0\ts,\al_1,\al_2\ts,\ts\ldots$
have been given in \cite[Section 4]{KN1}.
As well as in
\cite{KN1}, in the present article for any
$l=0\com1\com2\com\ts\ldots$ we
use the homomorphism $\Y(\gl_n)\to\Cr_{\ts l}$
$$
\be_{\ts l}\ts=\ts\al_{\ts l}\,\ts\om_n\ts=\ts
\pi_{n+l}\,\ts
\om_{n+l}\,\ts
\io_{\ts l}\,\ts
\om_n\,\ts
\tau_{\ts l}\ts.
$$
The last equality follows from the definition \eqref{al}
and the relation
$
\tau_{-\ts l}\,\ts\om_n
\ts=\ts
\om_n\,\ts\tau_{\ts l}\ts,
$
see \eqref{tauz} and \eqref{1.51}. The image
of any series \eqref{tser} under the homomorphism $\be_l$ can be expressed
in terms of quasideterminants \cite[Lemma 4.2]{BK} or
quantum minors \cite[Lemma 8.5]{BK}\ts; see also \cite[Lemma 1.5]{NT2}.
The reason for considering here the homomorphism $\be_{\ts l}$ rather than
$\al_{\ts l}$ will become clear when we state Proposition~\ref{arol},
see also \cite[Proposition 2.5]{C1}.

Consider the Lie algebra $\gl_m$ and its Cartan subalgebra $\h\ts$.
A weight $\mu\in\h^\ast$ is called \textit{polynomial\/} if
its labels $\mu_1\lcd\mu_m$ are non-negative integers such that
$\mu_1\ge\ldots\ge\mu_m\ts$. 
The weight $\mu\in\h^\ast$ is polynomial if and only if
for some non-negative integer $N$ the vector space
$$
\Hom_{\,\gl_m}(\ts L_{\ts\mu}\ts,(\CC^{\ts m})^{\ot N}\ts)\neq\{0\}\ts.
$$
Then 
$$
N=\mu_1+\ldots+\mu_m\ts.
\vspace{2pt}
$$
The irreducible $\gl_m$-module $L_{\ts\mu}$ of highest weight $\mu$
is then called a \textit{polynomial\/} module.
Then by setting $\mu_{m+1}=\mu_{m+2}=\ldots=0$
we get a partition 
$(\ts\mu_1,\mu_2\ts,\ts\ldots\ts)$
of $N$. When there is no confusion
with the polynomial weight of $\gl_m\ts$, this partition will be
denoted by $\mu$ as well. The maximal index $a$ with $\mu_a>0$
is then called the \textit{length\/} of the partition, and
is denoted by $\ell(\mu)\ts$. Note that here $\ell(\mu)\le m\ts$.
Further, let $\mup=(\mup_1,\mup_2\ts,\ts\ldots\ts)$ be 
the partition \textit{conjugate\/} to the partition $\mu\ts$.
By definition, here $\mup_b$ is equal to the maximal index $a$ such that
$\mu_a\ge b\ts$. Then $\mup_1=\ell(\mu)\ts$.

Let $\la$ and $\mu$ two polynomial weights of $\gl_m$ such that 
\begin{equation}
\label{lcon}
\ell(\lap)\le n+l
\quad\textrm{and}\quad
\ell(\mup)\le l\ts.
\end{equation}
Using the respective partitions, then
$\lap$ and $\mup$ can also be regarded as polynomial weights
of the Lie algebras $\gl_{n+l}$ and $\gl_{\ts l}$ respectively.
Denote by $L_{\ts\la}^{\ts\prime}$ and $L_{\ts\mu}^{\ts\prime}$
the corresponding irreducible highest weight modules
over $\gl_{n+l}$ and $\gl_{\ts l}\ts$.

Using the action of the Lie algebra $\gl_{\ts l}$ on 
$L_{\ts\la}^{\ts\prime}$ via its
embedding into $\gl_{n+l}$
as the second direct summand of the subalgebra 
$\gl_n\op\gl_{\ts l}\subset\gl_{n+l}$ take the vector space
\begin{equation}
\label{hll}
\Hom_{\,\gl_{\ts l}}(\ts L_{\ts\mu}^{\ts\prime}\ts,
L_{\ts\la}^{\ts\prime})\ts.
\end{equation}
The subalgebra $\Cr_{\ts l}\subset\U(\gl_{n+l})$ acts on this vector
space via the action of $\U(\gl_{n+l})$ on $L_{\ts\la}^{\ts\prime}\ts$. 
Moreover, the action of $\Cr_{\ts l}$ on \eqref{hll} is irreducible
\cite[Theorem 9.1.12]{D}. Hence the following identifications
of $\Cr_{\ts l\ts}$-modules are unique up to rescaling of the
vector \text{spaces\ts:}
$$
\Hom_{\,\gl_{\ts l}}(\ts L_{\ts\mu}^{\ts\prime}\ts,
L_{\ts\la}^{\ts\prime})\ts=
\vspace{4pt}
$$
$$
\Hom_{\,\gl_{\ts l}}(\ts L_{\ts\mu}^{\ts\prime}\ts,
\Hom_{\,\gl_m}(\ts L_{\ts\la}\ts,
\S\ts(\ts\CC^{\ts m}\ot\CC^{\ts n+l}\ts)))\ts=
\vspace{4pt}
$$
$$
\Hom_{\,\gl_{\ts l}}(\ts L_{\ts\mu}^{\ts\prime}\ts,
\Hom_{\,\gl_m}(\ts L_{\ts\la}\ts,
\S\ts(\ts\CC^{\ts m}\ot\CC^{\ts l}\ts)\ot
\S\ts(\ts\CC^{\ts m}\ot\CC^{\ts n}\ts)))\ts=
\vspace{2pt}
$$
\begin{equation}
\label{hlls}
\Hom_{\,\gl_m}(\ts L_{\ts\la}\ts,
L_{\ts\mu}\ot\S\ts(\ts\CC^{\ts m}\ot\CC^{\ts n}\ts))\ts.
\vspace{4pt}
\end{equation}
We used the classical identifications 
of modules over the Lie algebras $\gl_{n+l}$ and $\gl_m\ts$,
$$
L_{\ts\la}^{\ts\prime}=
\Hom_{\,\gl_m}(\ts L_{\ts\la}\ts,
\S\ts(\ts\CC^{\ts m}\ot\CC^{\ts n+l}\ts))
$$
and
$$
\Hom_{\,\gl_{\ts l}}(\ts L_{\ts\mu}^{\ts\prime}\ts,
\S\ts(\ts\CC^{\ts m}\ot\CC^{\ts l}\ts))=
L_{\ts\mu}
\vspace{4pt}
$$
respectively, see for instance \cite[Section 4.1]{H}.
We also use the decomposition
$$
\S\ts(\ts\CC^{\ts m}\ot\CC^{\ts n+l}\ts)=
\S\ts(\ts\CC^{\ts m}\ot\CC^{\ts l}\ts)\ot
\S\ts(\ts\CC^{\ts m}\ot\CC^{\ts n}\ts)\ts.
$$

By pulling back through the homomorphism $\be_l:\Y(\gl_n)\to\Cr_{\ts l}\ts$,
the vector space \eqref{hlls} becomes a module over the Yangian
$\Y(\gl_n)\ts$. On the other hand, the target vector
space $L_{\ts\mu}\ot\S\ts(\ts\CC^{\ts m}\ot\CC^{\ts n}\ts)$ in \eqref{hlls}
coincides with the vector space of the bimodule $\A_{\ts m}(L_{\ts\mu})$ over
$\gl_m$ and $\Y(\gl_n)$, see the beginning of Section 2. 
Using this bimodule structure, the vector space \eqref{hlls} becomes
another module over $\Y(\gl_n)\ts$.
But the next proposition shows that these two
$\Y(\gl_n)$-modules are the same. 
This proposition also makes \cite[Remark 12]{A}~more~precise.
We will give a direct proof of this proposition; 
another proof can be obtained by using \cite[Lemma 4.2]{BK}. 

\begin{proposition}
\label{arol}
The action of the algebra\/ $\Y(\gl_n)$ on the vector space \eqref{hlls}
via the homomorphism\/ $\be_l$ coincides with the action inherited from
the bimodule\/ $\A_{\ts m}(L_{\ts\mu})\ts$.
\end{proposition}

\begin{proof}
Consider the action of subalgebra $\Cr_{\ts l}\subset\U(\gl_{n+l})$
on $\S\ts(\ts\CC^{\ts m}\ot\CC^{\ts n+l}\ts)\ts$.
The Yangian $\Y(\gl_n)$ acts on this vector space via the homomorphism
$\be_l:\Y(\gl_n)\to\Cr_{\ts l}\ts$. Identify this vector space with 
$\P\ts(\ts\CC^{\ts m}\ot\CC^{\ts n+l}\ts)$ as before.
Using the decomposition
\begin{equation}
\label{pmnl}
\P\ts(\ts\CC^{\ts m}\ot\CC^{\ts n+l}\ts)=
\P\ts(\ts\CC^{\ts m}\ot\CC^{\ts l}\ts)\ot
\P\ts(\ts\CC^{\ts m}\ot\CC^{\ts n}\ts)
\end{equation}
we will show that
for any $s=0,1,2,\ts\ldots$ and $i\com j=1\lcd n$ the generator
$T_{ij}^{\ts(s+1)}$ of $\Y(\gl_n)$ then acts on the vector space \eqref{pmnl} 
as the element \eqref{combact} of the algebra 
$\U(\gl_m)\ot\PD\ts(\CC^{\ts m}\ot\CC^{\ts n})\ts$.
Proposition \ref{arol} will thus follow from Proposition \ref{dast}.

For any $i\com j=1\lcd n+l$
the element $E_{ij}\in\U(\gl_{n+l})$ acts on
$\P\ts(\ts\CC^{\ts m}\ot\CC^{\ts n+l}\ts)$ as the differential operator
$$
\sum_{c=1}^m\,
x_{ci}\,\d_{\ts cj}\,.
$$
Consider the $(n+l)\times(n+l)$ matrix whose $ij$-entry is
$$
\de_{ij}+(u-l)^{-1}\,
\sum_{c=1}^m\,
x_{ci}\,\d_{\ts cj}\,.
$$
Write this matrix and its inverse as the block matrices
$$
\begin{bmatrix}
\,A&B\,\\\,C&D\,
\end{bmatrix}
\quad\textrm{and}\quad
\begin{bmatrix}\ \At&\Bt\ \\\ \Ct&\Dt\
\end{bmatrix}
$$
where the blocks $A,B,C,D$ and $\At,\Bt,\Ct,\Dt$ are matrices of sizes
$n\times n$, $n\times l$, $l\times n$, $l\times l$ respectively. 
The action of the algebra $\Y(\gl_n)$ on the vector space
$\P(\ts\CC^{\ts m}\ot\CC^{\ts n+l}\ts)$ via the homomorphism
$\be_l:\Y(\gl_n)\to\Cr_{\ts l}$ can now be described by assigning
to the series $T_{ij}(u)$ with $i\com j=1\lcd n$ the 
$ij$-entry of the matrix $\At^{\,-1}$.

Consider the $(n+l)\times m$ matrix whose $ic\ts$-entry is 
the operator of multiplication by $x_{ci}$ in 
$\P(\ts\CC^{\ts m}\ot\CC^{\ts n+l}\ts)$ on the left. Write this matrix as
$$
\begin{bmatrix}
\,P\,
\\
\,\Pb\,
\end{bmatrix}
$$
where the blocks $P$ and $\Pb$ are matrices of sizes
$n\times m$ and $l\times m$ respectively.
Consider the $m\times(n+l)$ matrix whose $cj$-entry is 
the operator $\d_{cj}\ts$. Write it as
$$
\begin{bmatrix}
\,Q\,\,\Qb\,\ts
\end{bmatrix}
$$
where the blocks $Q$ and $\Qb$ are matrices of sizes
$m\times n$ and $m\times l$ respectively. Then
$$
\begin{bmatrix}
\,A&B\,\\\,C&D\,
\end{bmatrix}
\,=\,
1+(u-l)^{-1}
\begin{bmatrix}
\,P\,
\\
\,\Pb\,
\end{bmatrix}
\begin{bmatrix}
\,Q\,\,\Qb\,
\end{bmatrix}
\,=\,
1+(u-l)^{-1}
\begin{bmatrix}\,
\ P\ts Q\,&\,P\ts\Qb\ 
\\ 
\ \Pb\ts Q\,&\,\Pb\ts \Qb\  
\end{bmatrix}
$$ 
where $1$ stands for the $(n+l)\times(n+l)$ identity matrix.
Using Lemma \ref{lemma1}, we get
$$
\At^{\,-1}=\,A-B\ts D^{-1}\ts C
\,=\,
\vspace{4pt}
$$
$$
1+(u-l)^{-1}\ts P\ts Q-
(u-l)^{-2}\ts 
P\ts\Qb\,\bigl(\ts1+(u-l)^{-1}\ts\Pb\ts\Qb\,\bigr)^{-1}\ts\Pb\ts Q
\,=\,
$$
\begin{equation}
\label{pulq}
1+P\ts(\ts u-l+\Qb\ts\Pb\,\bigr)^{-1}\ts Q
\vspace{2pt}
\end{equation}
where $1$ now stands for the $n\times n$ identity matrix. 

Now consider the $m\times m$ matrix $u-l+\Qb\ts\Pb$ appearing in
\eqref{pulq}. Its $ab$-entry is
$$
\de_{ab}\ts(u-l)\ts+\sum_{k=1}^l\,
\d_{\ts a,n+k}\,x_{\ts b,n+k}
\,=\,
\de_{ab}\ts u\ts-\sum_{k=1}^l\,
x_{\ts b,n+k}\,\d_{\ts a,n+k}\,.
$$
Observe that the last displayed sum over $k=1\lcd l$ corresponds 
to the action of the element $E_{\ts ba}\in\U(\gl_m)$ 
on the first tensor factor in the decomposition \eqref{pmnl}. 
Denote by $Y_{ab}(u)$ the $ab$-entry of the matrix inverse to 
$u-l+\Qb\ts\Pb$. The $ij$-entry of the 
matrix \eqref{pulq} can then be written as the sum
$$
\de_{ij}\ts+\sum_{a,b=1}^m x_{ai}\,Y_{ab}(u)\,\d_{\ts bj}
\,=\,
\de_{ij}\ts+\sum_{a,b=1}^m Y_{ab}(u)\,x_{ai}\,\d_{\ts bj}\,.
$$
Using the observation made in the end of Section 1 now completes the proof.
\end{proof}

Note that our proof of Proposition \ref{arol} remains
valid in the case $l=0$. In this case we assume that
$\gl_{\ts l}=\{0\}\ts$. For any positive integer $l$ 
consider the homomorphism
$\U(\gl_m)\to\PD\ts(\CC^{\ts m}\ot\CC^{\ts l}\ts)$
corresponding to the action of $\gl_m$  
on the first tensor factor in the decomposition \eqref{pmnl}.
The kernels of all these homomorphisms 
for $l=1\com2\com\ts\ldots$
have the zero intersection.
Thus, independently of Proposition~\ref{dast},
our proof of Proposition~\ref{arol} shows that 
for any positive integer $m$ a homomorphism of associative algebras
\eqref{lastrem} can be defined
by mapping $T_{ij}^{\ts(s+1)}$ to the sum \eqref{combact}.

Further, for any given polynomial weights $\la$ and $\mu$ of
$\gl_m$ we can choose an integer $l$ large
enough to satisfy \eqref{lcon}. Then
the algebra $\Cr_{\ts l}$ acts on the vector
space \eqref{hlls} irreducibly, while the central elements of
$\U(\gl_{n+l})$ act on \eqref{hlls} via multiplication by scalars.
So Proposition \ref{arol} has a corollary;
cf.\ \cite[Theorem~10]{A}.

\begin{corollary}
The action of the algebra\/ $\Y(\gl_n)$ on the vector space \eqref{hlls}
inherited from the bimodule\/ $\A_{\ts m}(L_{\ts\mu})$ is irreducible
for any polynomial weights $\la$ and $\mu$ of\/ $\gl_m\ts$. 
\end{corollary}

It is well known 
that the vector space \eqref{hll}
is not zero if and only if 
\begin{equation}
\label{llcon}
0\le\la_a-\mu_a\le n
\quad\textrm{for every}\quad
a=1\lcd m\ts;
\end{equation} 
Hence the space \eqref{hlls}
is also not zero if and only if the inequalities \eqref{llcon} hold.
For further details on the irreducible representations of the Yangian 
$\Y(\gl_n)$ of the form \eqref{hll} see for instance \cite[Section 4]{M} and
\cite[Section 2]{NT2}.



\begin{thebibliography}{[MNO]}

\bibitem[A]{A}
{T.\,Arakawa},
\textit{Drinfeld functor and finite-dimensional representations of the 
Yangian},
{Commun. Math. Phys.}
\textbf{205}
(1999), 
1--18.

\bibitem[AS]{AS}
{T.\,Arakawa and T.\,Suzuki},
\textit{Lie algebras and degenerate affine Hecke algebras of type $A$},
{J. Algebra}
\textbf{209}
(1998),
288--304.

\bibitem[AST]{AST}
{T.\,Arakawa, T.\,Suzuki and A.\,Tsuchiya},
\textit{Degenerate double affine Hecke algebras and conformal field theory},
{Progress Math.}
\textbf{160}
(1998),
1--34.

\bibitem[BK]{BK}
{J.\,Brundan and A.\,Kleshchev},
\textit{Parabolic presentations of the Yangian $Y(\gl_n)$},
{Commun. Math. Phys.}
\textbf{254}
(2005),
191--220.

\bibitem[C1]{C1}
{I.\,Cherednik},
\textit{ A new interpretation of Gelfand-Zetlin bases},
{Duke Math. J.}
{\bf 54}
(1987),
563--577.

\bibitem[C2]{C2}
{I.\,Cherednik},
\textit{Lectures on Knizhnik-Zamolodchikov
equations and Hecke algebras},
{Math. Soc. Japan Memoirs}
{\bf 1}
(1998),
1--96.

\bibitem[D]{D}
J.\,Dixmier,
\textit{Alg\`ebres enveloppantes},
Gauthier-Villars, Paris, 1974.

\bibitem[D1]{D1}
{V.\,Drinfeld},
\textit{Hopf algebras and the quantum Yang-Baxter equation},
{Soviet Math.\,Dokl.}
\textbf{32}
(1985),
254--258.

\bibitem[D2]{D2}
{V.\,Drinfeld},
\textit{Degenerate affine Hecke algebras and Yangians},
{Funct. Anal. Appl.}
\textbf{20}
(1986),
56--58.

\bibitem[H1]{H1}  
{R.\,Howe},
\textit{Remarks on classical invariant theory},
Trans. Amer. Math. Soc.
\textbf{313}  
(1989), 
539--570. 

\bibitem[H2]{H}  
{R.\,Howe},
\textit{Perspectives on invariant theory: 
Schur duality, multiplicity-free actions and beyond},
Israel Math. Conf. Proc.
\textbf{8}
(1995),
1--182.

\bibitem[K]{K}
{S.\,Khoroshkin},
\textit{Extremal projector and dynamical twist},
{Theoret. Math. Phys.}
\textbf{139} 
(2004),
582--597.

\bibitem[KN]{KN1}
{S.\,Khoroshkin and M.\,Nazarov},
\textit{Yangians and Mickelsson algebras I\/},
Transformation Groups 
\textbf{11} 
(2006),
to appear.

\bibitem[KO]{KO}
{S.\,Khoroshkin and O.\,Ogievetsky},
\textit{Mickelsson algebras and Zhelobenko operators},
\texttt{arXiv:math.QA/0606259}.

\bibitem[L]{L}  
G.\,Lusztig,
\textit{Affine Hecke algebras and their graded version},
{J. Amer. Math.Soc.}
{\bf2} 
(1989),   
599--635.

\bibitem[M]{M}
{A.\,Molev},
\textit{Yangians and transvector algebras},
Discrete Math.
\textbf{246}
(2002),
231--253.

\bibitem[M1]{M1}
{J.\,Mickelsson},
\textit{Step algebras of semi-simple subalgebras of Lie algebras},
Reports Math. Phys. 
\textbf{4}
(1973),
307--318.

\bibitem[M2]{M2}
{J.\,Mickelsson},
\textit{On irreducible modules of a Lie algebra which are composed of
finite-dimensional modules of a subalgebra},
{Ann. Acad. Sci. Fenn. Ser. A I}
\textbf{598}
(1975),
1--16.

\bibitem[MNO]{MNO}
{A.\,Molev, M.\,Nazarov and G.\,Olshanski},
\textit{Yangians and classical Lie algebras},
{Russian Math. Surveys}
\textbf{51} 
(1996),
205--282.

\bibitem[N]{N1}
{M.\,Nazarov},
\textit{Yangian of the queer Lie superalgebra},
{Commun. Math. Phys.}
\textbf{208}
(1999), 
195--223.

\bibitem[NT1]{NT2}
{M.\,Nazarov and V.\,Tarasov},
\textit{Representations of Yangians with Gelfand-Zetlin bases},
{J. Reine Angew. Math.}
\textbf{496}
(1998),
181--212.

\bibitem[NT2]{NT1}
{M.\,Nazarov and V.\,Tarasov},
\textit{On irreducibility of tensor products of Yangian modules},
{Internat. Math. Res. Notices}
(1998),
125--150.

\bibitem[O1]{O1}
{G.\,Olshanski},
\textit{Extension of the algebra $U(g)$ for infinite-dimensional classical
Lie algebras $g$, and the Yangians $Y(gl(m))$},
{Soviet Math. Dokl.}
\textbf{36}
(1988),
569--573.

\bibitem[O2]{O2}
{G.\,Olshanski},
\textit{Representations of infinite-dimensional classical groups,
limits of enveloping algebras, and Yangians},
{Adv. Soviet Math.}
\textbf{2}
(1991),
1--66.

\bibitem[ST]{ST}
Y.\,Smirnov and V.\,Tolstoy,
\textit{Extremal projectors for usual, super and quantum algebras
and their use for solving Yang-Baxter problem}, 
Selected Topics in Mathematical Physics, 
World Scientific, Teaneck, 1990, pp. 347--359. 

\bibitem[TV]{TV2}
{V.\,Tarasov and A.\,Varchenko},
\textit{Duality for Knizhnik-Zamolodchikov and dynamical equations},
{Acta Appl. Math.}
{\bf 73}
(2002), 
141--154. 

\bibitem[Z]{Z}
D.\,Zhelobenko, 
\textit{Extremal cocycles on Weyl groups},
{Funct. Anal. Appl.}
{\bf 21}
(1987),
183--192.

\end{thebibliography}
\end{document}